\documentclass[acmsmall,screen,nonacm]{acmart}

\AtBeginDocument{%
  \providecommand\BibTeX{{%
    \normalfont B\kern-0.5em{\scshape i\kern-0.25em b}\kern-0.8em\TeX}}}

\copyrightyear{2022}
\acmYear{2022}
\setcopyright{rightsretained}
\acmJournal{TSAS}
\acmYear{2022} \acmVolume{00} \acmNumber{JA} \acmArticle{00} \acmMonth{6} \acmPrice{}\acmDOI{10.1145/3544779}

\newcommand{\algo}{{\textsf{SPATIAL}}}
\newcommand{\spatial}{{\textsf{SOP}}}

\def\modelA{{\textsf{BAA}}}
\def\modelB{{\textsf{BCAA}}}
\def\modelC{{\textsf{AIO}}}

\def\lcps{{\bf X}}
\def\fcps{{\bf Y}}

\def\graph{{\mathcal{G}}}
\def\vertexset{{\mathcal{V}}}
\def\edgeset{{\mathcal{E}}}
\def\indexset{{\mathcal{I}}}

\def\adjacency{{\mathbf{W}}}
\def\assignment{{\mathbf{X}}}

\def\flip{\textsf{Flip}}
\def\adjacent{W}
\def\assign{{X}}
\def\terminal{{T}}

\def\capacity{{\mu}}

\def\objfunc{\mathbf{\mathcal{J}}}

\def\u{{\vphantom g}_}
\def\level{\mathrm{L}}
\def\nX{{\u \level}\!p}

\usepackage{wrapfig}
\usepackage{nicefrac}
\usepackage{multirow}
\usepackage[boxed, ruled, vlined]{algorithm2e}
\usepackage{algpseudocode}
\usepackage{subcaption}
\usepackage[noabbrev,capitalise]{cleveref}

\begin{document}


\title{
Memetic algorithms for {S}patial {P}artitioning problems
}


\author{Subhodip Biswas}
\authornote{This work is a part of the author's doctoral dissertation \cite{biswas2022phd}.
\newline
This research was supported in parts by National Science Foundation (NSF) grants DGE-1545362 and IIS-1633363. 
}
\email{subhodip@vt.edu}
\orcid{0000-0002-7065-2966}
\author{Fanglan Chen}
\orcid{0000-0001-7610-3178}
\email{fanglanc@vt.edu}
\affiliation{%
  \institution{Virginia Tech}
  \state{Virginia}
  \country{USA}
}

\author{Zhiqian Chen}
\orcid{0000-0003-4112-9647}
\affiliation{%
  \institution{Mississippi State University}
  \city{Starkville}
  \state{Mississippi}
  \country{USA}}
\email{zchen@cse.msstate.edu}

\author{Chang-Tien Lu}
\orcid{0000-0003-3675-0199}
\email{clu@vt.edu}
\author{Naren Ramakrishnan}
\orcid{0000-0002-1821-9743}
\email{naren@cs.vt.edu}
\affiliation{%
  \institution{Virginia Tech}
  \state{Virginia}
  \country{USA}
}

\renewcommand{\shortauthors}{Biswas et al.}

\begin{abstract}
Spatial optimization problems (SOPs) are characterized by spatial relationships governing the decision variables, objectives, and/or constraint functions. In this article, we focus on a specific type of SOP called spatial partitioning, which is a combinatorial problem due to the presence of discrete spatial units. Exact optimization methods do not scale with the size of the problem, especially within practicable time limits. This motivated us to develop population-based metaheuristics for solving such SOPs. However, the search operators employed by these population-based methods are mostly designed for real-parameter continuous optimization problems. For adapting these methods to SOPs, we apply domain knowledge in designing spatially-aware search operators for efficiently searching through the discrete search space while preserving the spatial constraints. To this end, we put forward a simple yet effective algorithm called \textsf{s}warm-based s\textsf{pa}tial meme\textsf{ti}c \textsf{al}gorithm (\algo) and test it on the school (re)districting problem. Detailed experimental investigations are performed on real-world datasets to evaluate the performance of \algo. Besides, ablation studies are performed to understand the role of the individual components of \algo. Additionally, we discuss how \algo~is helpful in the real-life planning process and its applicability to different scenarios and motivate future research directions.
\end{abstract}

\maketitle

\section{Introduction}
Spatial optimization has been an active research area, especially in disciplines such as economics, engineering, environmental studies, geography, operational research, and regional science. \citet{churchSpatial} noted that ``spatial optimization involves identifying how land use and other activities should be arranged and organized spatially in order to optimize efficiency or some other quality measure.'' It includes many districting, layout, location and network problems that involve design, operations, and planning~\cite{cao20171}. Solving a spatial optimization problem (\spatial) is equivalent to searching for an optimal assignment of a set of discrete spatial units representing some geographic areas such that some well-defined objectives and/or constraints are satisfied. Alternatively, we can also define this as the use of mathematical or computational techniques for finding solutions to geographic decision problems subjected to design constraints~\cite{spatialopt}. The optimization variables in a \spatial~relate to the decision being made with the objective function quantifying the quality of the decision. The constraints impose a set of necessary design considerations that needs to be satisfied. The functions and constraints usually encode spatial properties/topological relationships, including adjacency, contiguity (connectivity), similarity (distance), shape (compactness), and so on~\cite{tong2012spatial}.

Broadly, \spatial s~are usually classified as either a \textit{selection} or a \textit{partitioning} problem~\cite{xiao2008unified}. 
Spatial selection problems identify a subgroup of spatial units. Additional spatial constraints need to be satisfied for certain problems, while others only impose continuations on the selected spatial units. Spatial partitioning problems, on the other hand, seek to group the spatial units into a number of districts or territories. For instance, consider the districting problem~\cite{kalcsics2019districting}, where the objective is to partition a geographic area into groups of contiguous districts (regions) such that each district is balanced with respect to some activity measure, like residing population. 
Due to the discrete nature of spatial units, \spatial s suffer from combinatorial explosion, i.e., the phenomenon where the computing time cost to find the optimal solution of a \textsf{NP}-hard problem increases exponentially with the problem size~\cite{megiddo1984complexity,gilbert1985multiobjective}. Thus, exact optimization techniques like Integer Programming (IP) or mixed-integer Linear Programming (MILP) cannot solve the problem optimally under practical time constraints~\cite{glover2015metaheuristics}. This is why researchers often resort to using approximation methods like \textit{heuristics} and \textit{metaheuristics} since these methods can find good, but not necessarily optimal, solutions to the problem in a reasonable time. Thus, computational efficiency is the key to designing these methods for solving SOPs~\cite{xiao2015gis}.

Heuristic methods are mainly designed for solving a particular problem. For instance, let us consider the $p-$Median Problem ($p$MP) in location sciences, where the aim is to find facilities or services on $p$ nodes of a network such that the distance from each node to its nearest facility or service node is minimized~\cite{hakimi1964optimum}. \citet{teitz1968heuristic} proposed a heuristic called the vertex-exchange algorithm to solve this problem. This heuristic starts with a randomly selected subset of $p$ nodes from the network and keeps flipping these nodes with unselected ones until such exchange can no longer improve the quality of the solution.
In contrast to heuristics, metaheuristics refer to a general problem-solving framework that is composed of a set of high-level problem-independent instructions or strategies for developing heuristic optimization algorithms~\cite{sorensen2013metaheuristics}. Some well-known examples of metaheuristics include evolutionary algorithms~\cite{holland1992adaptation}, simulated annealing~\cite{kirkpatrick1983optimization}, tabu search~\cite{glover1998tabu}, variable neighborhood search~\cite{vns}, etc. Oftentimes, these methods are inspired by some natural processes and they can be adapted to solve different kinds of problems. Hence, metaheuristics have become a popular choice amongst practitioners and researchers for solving medium to large instances of combinatorial optimization problems~\cite{blum2011hybrid}, specially in location sciences~\cite{mladenovic2007p}.

Motivated by this, we devise a simple, easy-to-use population-based metaheuristic inspired by the emerging field of Swarm Intelligence.\footnote{Swarm Intelligence is defined as the study and design of computational optimization techniques based on the collective intelligence emerging from a large population of search agents with simple behavioral patterns for communication and interaction. These methods  instantiate search moves the closely mimic the complex social behavior of animals such as ant colonies, beehives, bird flocks, and so on \cite{del2019bio}.}
In particular, our approach builds on top of the Artificial Bee Colony (ABC) algorithm that is based on the foraging behavior of the swarm of bees~\cite{akay2021survey}. It maintains multiple trial solutions to the \spatial~under consideration and combines a local search technique with a spatially-aware recombination operator resulting in what is commonly known as a \textit{memetic algorithm}~\cite{moscato2010modern}. Hence the name \textbf{s}warm-based s\textbf{pa}atial meme\textbf{ti}c \textbf{al}gorithm (\algo). The search move of the algorithm is modified to explore the discrete search space while preserving the spatial relationships/constraints.  preliminary version of this work appeared in \cite{spatial}. In this paper, we present further additions based on theoretical and experimental investigations of our framework. The summary of the extensions and contributions are elucidated below.

\begin{itemize}
    \item An overview of background details is provided in \cref{sec:background}, specially in context of spatial partitioning problems like districting in \cref{sec:sop}. We then show in \Cref{sec:graph} spatial partitioning problems are accompanied by an underlying graph structure which enables the problem to be solved as a graph partitioning problem. \Cref{sec:meta} briefly reviews graph-partitioning techniques that motivate our algorithmic approach and the role of domain knowledge in algorithm design.
    
    \item \Cref{sec:optimization} defines the optimization problem corresponding to a generalized spatial partitioning by leveraging the notions from graph partitioning. We also show how the given formulation can be adapted to problems like school districting in \Cref{sec:school} and provide some pointers on how to adapt other types of spatial partitioning problems using the given framework.
    
    \item A detailed outline of the \algo~method is proposed in \Cref{sec:method}. In particular, we focus on the two improvement steps: \textit{local search} and \textit{spatially-aware recombination}. Additionally, the relationship between locals search and the sampling of partitions based on the theory of Markov Chain is discussed.
    This is followed by an in-depth discussion on how the recombination operation efficiently searches for solutions in the discrete search space.
    
    \item We used the dataset for the school year 2020-21 here as compared to 2019-20 used in the previous work \cite{spatial}. At the onset of the pandemic in 2020, many parents unenrolled their children from public schools in these two districts thereby creating a serious imbalance between the student population and the school capacity. This presented more challenging problem instances to work with.

    \item We included more sophisticated baseline methods and performed an exhaustive comparison in \Cref{sec:exp}. Additional discussions on how solution initialization affects performance and how the algorithms can be made more efficient using alternative measures. We also include a case study showcasing the applicability of \algo~in real-life planning.
\end{itemize}


\section{Background}
\label{sec:background}
This section provides a basic outline of the background details necessary for understanding this research. Firstly, an introduction to spatial partitioning problems are provided in~\cref{sec:sop}. This is followed by \Cref{sec:graph} highlighting how graphs can encode relationships between the spatial units in most of the spatial problems. In fact, the graph-based representation can be used to pose spatial partitioning as a graph partitioning problem. A brief review of different graph partitioning approaches is provided in \Cref{sec:meta}.
 
\subsection{Spatial partitioning problems}
\label{sec:sop}
The field of spatial optimization is firmly rooted in the classical works on graph theory, where mathematical formalism and theories about spatial arrangement and movement were made.
In spite of its historic origin, the term \textit{spatial optimization} first appeared in the literature during the late 1960s and the early 1970s, when a series of articles made use of the term within the context of regional science~\cite{mathieson1969soviet, alao1971two, marble1972landuse}. Interestingly, this research domain has followed the developments in Computer Science, which is rich in works on graphs and other discrete data structures. 

The term ``spatial optimization'' was initially used by~\citet{ghosh1984location} to describe a set of location-allocation problems, namely the warehouse-location problem~\cite{manne1964plant}, the $p-$MP~\cite{revelle1970central}, and the location set-covering problems~\cite{toregas1972optimal}. Similarly, researchers in resource management use the term spatial optimization to refer to optimization models for allocating various land-use areas within a forest~\cite{hof1992spatial,hof1998spatial}.
In fact, spatial optimization problems appear in different disciplines in different context$-$location sciences~\cite{hale2003location, yao2018spatial, laporte2019introduction}, regionalization~\cite{lankford1969regionalization, surveysupervised}, spatial data mining~\cite{han2001spatial, mennis2009spatial}, territory design~\cite{zoltners1983sales, kalcsics2005towards, salazar2011new}, etc. 
Most of these developments have a commonality$-$optimizing an objective function subjected to a number of constraints (spatial or aspatial) that define the feasibility of solutions. These include a large number of simplistic variations of well-known problems, including the location set covering, $p-$median, simple plant location problems, etc. Such normative location problems have garnered a lot of interest in spatial optimization, especially from the research community and the industry. In fact, \spatial s are too broad to be addressed all at one time and is outside the scope of this work. In this article, we focus on a spatial partitioning problem called \textit{territory design problem}, often popularly known as \textit{districting} or \textit{zone design} in location sciences~\cite{laporte2019introduction}. Note that the term \textit{redistricting} is also used to refer to these problems. However, redistricting actually means rearrangement of existing territories. For the sake of clarity, we shall use the term ``districting'' to refer to these problems.

\subsubsection{Districting problems}
\label{sec:districting}
Districting is a sub-field of discrete optimization involving some form of partitioning decision. In a typical districting problem, a set of smaller geographic areas, called \textit{basic units} or \textit{spatial units}, are group together into larger geographic areas, called \textit{districts} or \textit{territories}, such that they meet a series of planning criteria and requirements as specified from the application or context~\cite{rios2020research}. Districting problems arise in different real-life application domains, ranging from political districting over the design of districts for police patrols, schools, social facilities, waste collection, or winter services, to sales and service territory design~\cite{kalcsics2019districting}. 

Based on the application domain, each category of districting problem is unique from the perspective of modeling, criteria, or constraints. Nonetheless, several common criteria, including balance, compactness, unique assignment, can be generally applied to most districting problems. ``Balance'' means that a total amount of resources need to be fairly allocated among the districts. The term resources implies a particular attribute or multiple attributes of each spatial unit. Examples include the number of customers, product demand, population size, workload, etc. ``Compactness'' aims to obtain districts composed of basic units with geographical proximity, which can be optimized by minimizing a dispersion function that measures how tightly the area of a district is packed within its perimeter. 
Unique assignment indicates that each spatial unit must be assigned to only one district, and this requirement assures a complete partitioning of all the spatial units. Additionally, territory ``contiguity'' needs to be considered while designing the districts/territories.

Interestingly, there is no single approach to model the aforementioned criteria. Therefore, existing works have studied various methods to represent these requirements. During the 1960s to 1980s, the majority of research has focused on sales territory design~\cite{zoltners1983sales} and political districting problems~\cite{ricca2011political}. In the last three decades, a lot of studies on other applications such as distribution territory design, service-related districting, and, more recently, districting in health care~\cite{yanik2020review}, have emerged.
Scant attention has been paid to the problem of school districting till now.

\paragraph{School districting}

\begin{figure}
    \centering
    \includegraphics[keepaspectratio, width=0.5\textwidth]{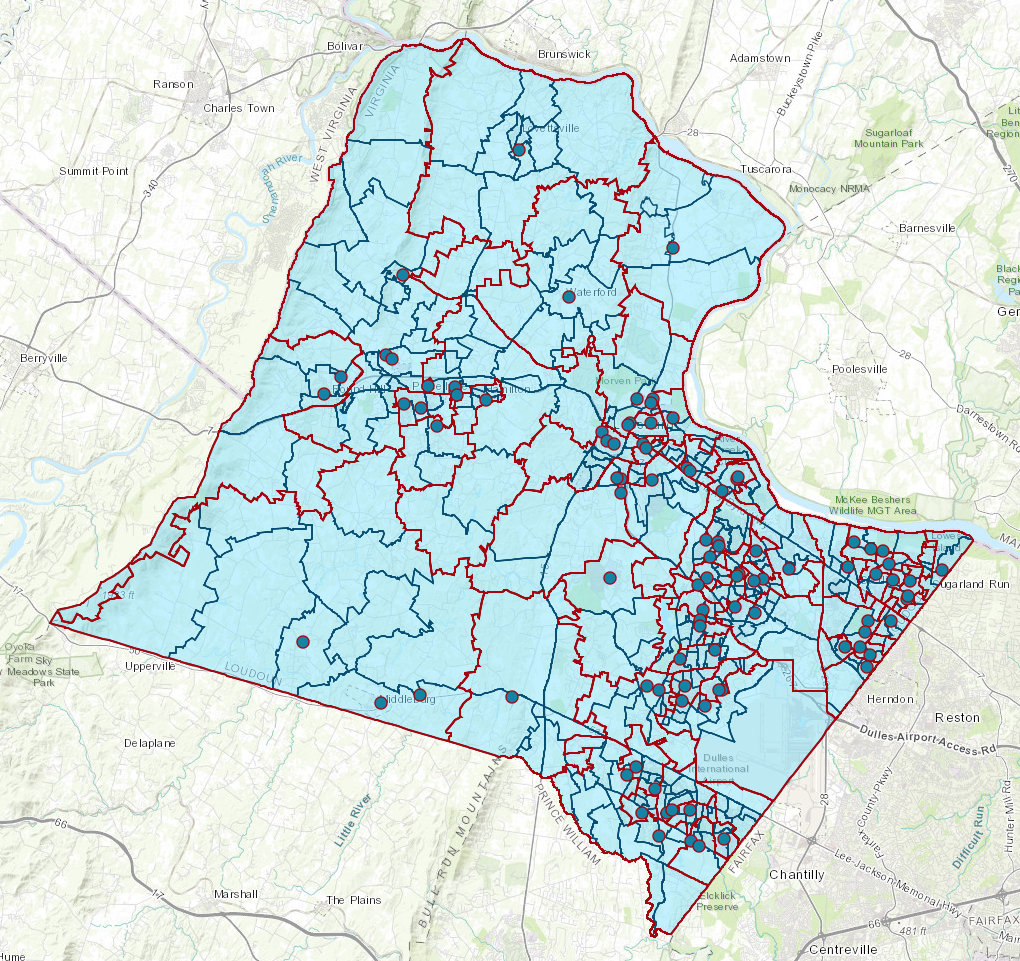}
    \caption{ArcGIS visualization of the school district of Loudoun county, VA, USA. The smaller polygons (with blue-colored border) represent the SPAs, while the larger polygons (with brown-colored border) represent the SAZs of elementary schools. The dark blue dots represent the locations of all the public schools.}
    \label{fig:schooldistrict}
\end{figure}

In countries like the US, \textit{school districts} play a vital role in the operation of the public school systems. A school district is an administrative unit for overseeing the jurisdiction of public schools and represent a large geographical region that is coterminous with the boundary of a county, city, or a subdivision. The spatial configuration of a school district is composed of smaller spatial units called planning zones or student planning areas (SPAs). These SPAs are grouped to form a larger geographically-contiguous area, called the \textit{school attendance zone} (SAZ), which defines the boundary of a school. The schools at every grade-level (elementary, middle, and high) have a well-defined boundary often arranged in a hierarchical manner. In a school district, the rule of thumb is that students attend the school assigned to their residing SPA. Figure~\ref{fig:schooldistrict} illustrates a map of a school district along the school locations, school boundaries and constituent SPAs. Note that in districting problems, a large geographical area, like a county, is partitioned in multiple ``districts'' or territories. However, in school districting, the term ``school district'' refers to the entire geographic area, like county or city. To avoid this confusion, one can imagine a school district being composed of contiguous districts or territories, each of which represent a school boundary.

School districting is the process by which the boundaries of public schools (within a school district) are adjusted/redrawn in response to projected growth/decline of student enrollment, change in school capacities, opening/closure of a school, etc. This is an annual/biannual event that involves the school boundary planners, board members, parents, and other stakeholders, and takes up a significant amount of time in reaching a consensus about the final districting plan to implement. Multiple factors (geographic, economic, social) are considered in deciding the school boundaries, thereby making school districting a technically and socially challenging to solve. The  complexity of the process piqued the curiosity of the research community, specially in Operations Research.

\citet{sutcliffe1984goal} summarized the work in this direction till the early 1980s. Since then not many works have been reported in this direction~\cite{schoepfle1989fast,ferland1990decision,lemberg2000school}. Among the few, \citet{schoepfle1989fast} introduced the term \textit{Generic School Districting Problem}, which refers to a class of school boundary problems involving allocation of students to schools while minimizing a cost/distance function, subject to a set of balancing/activity constraints.
Most of the approaches to school redistricting usually solves a continuous LP or a derived transportation problem in order to get an optimal or a near-optimal solution (which requires split-resolution)~\cite{belford1972network,liggett1973application}. However, the computational bottleneck does not allow these methods to scale.
\citet{caro2004school} proposed an IP approach to the school districting problem by minimizing the total distance travelled by the students. Their model was inspired by the sales territory alignment model of~\cite{zoltners1983sales} and was perhaps the first approach to account for all the problem-specific constraints, including connectivity. Their model was applied to only 2 (out of 22) clusters that the school district of the City of Philadelphia is divided into. 

\subsection{Graph-based representation in spatial optimization}
\label{sec:graph}

\begin{wrapfigure}{R}{0.5\textwidth}
  \begin{center}
    \includegraphics[keepaspectratio, width=0.48\textwidth]{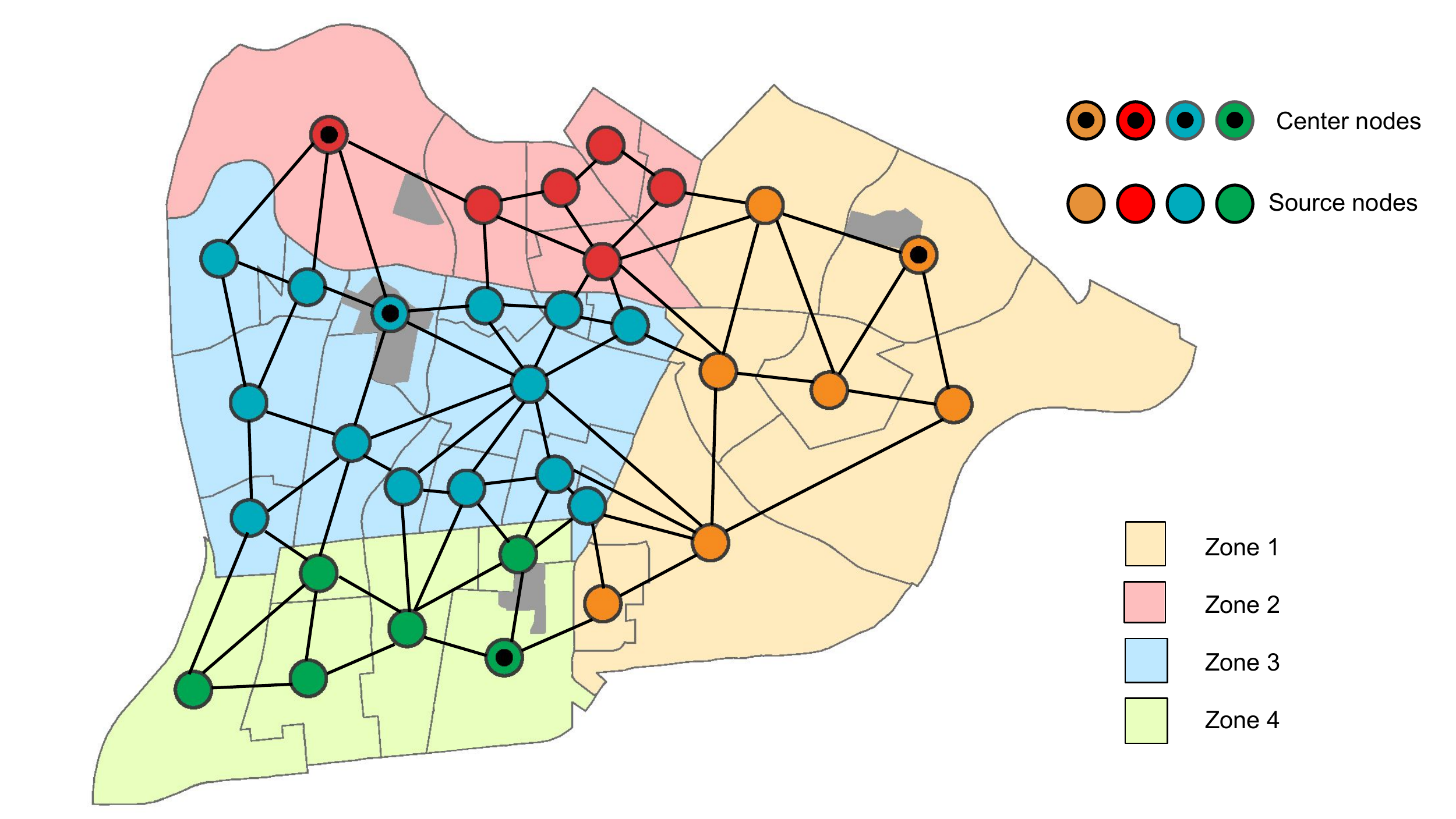}
    \caption{A large geographical area composed of smaller spatial units can be represented by a planar connected graph. The color coding represents an instance of spatial partitioning, where the each color correspond a territory/zone obtained by partitioning the graph into connected subgraphs.}
    \label{fig:graph}
    \end{center}
\end{wrapfigure}

A geographical area composed of smaller-sized spatial units can be represented as nodes of a graph $\graph = \bigl( \vertexset, \edgeset \bigr)$, where $\vertexset=\{v\u 1, v\u 2, \ldots, v\u N\}$ is the set of nodes representing the $N$ smaller-sized spatial units and $\edgeset$ is the set of edges connecting adjacent nodes. $\graph$ is commonly called the \textit{contiguity graph} or the \textit{dual graph}. It is a planar connected graph with the nodes encoding the spatial entities and the edges capturing the spatial adjacency relationship between the entities. 
An edge connects two nodes if their corresponding spatial units share a common boundary (more than a single point). \Cref{fig:graph} illustrates a toy example depicting the graph-based encoding of a geographical area.

A node $v$ may also be represented by its index, i.e., $v \in \{1, 2, \ldots, N\}$. These nodes may have features, i.e., $\mathbf{F} \left( \graph \right) = \left\{{F}_1, {F}_2, \ldots, {F}_N \right\}$, where ${F}\u{v}$ is the set of features corresponding to node $v$. Let ${F}\u{v}$ be represented by a tuple
$\bigl(L\u v, A\u{v} \bigr)$, where $L\u v=\left[(x_1,y_1), (x_2,y_2), \ldots, (x_t,y_t), (x_1,y_1)\right]$ is the list of geographic coordinates (latitude-longitude) that define the boundary polygon of the $v$\textsuperscript{th} spatial unit and $A\u{v}$ is the vector of feature values. 
Usually, a similarity matrix $\adjacency \left( \graph \right) = \bigl( \adjacent\u {uv} \bigr)\u {u=1,\dots,N\; v=1,\ldots,N}$ captures the relationship between any pair of nodes. Popular choices for the similarity metric include the distance function or the binary adjacency matrix. Similarity can also be defined for edges connecting adjacent nodes.

Solving a \spatial~involves the search for a feasible solution $\graph' = \bigl( \vertexset', \edgeset' \bigr)$, such that the spatial configurations of $\vertexset' \subseteq \vertexset$ and $\edgeset'\subseteq \edgeset$ satisfy predefined problem criteria/constraints while minimizing certain objective(s). In spatial partitioning problems like districting, we seek a $K-$partition of the graph $\graph$, i.e., $ \edgeset' \subset \edgeset$ and $\vertexset' = \vertexset$ such that the nodes in $\vertexset'$ are grouped into $K$ connected subgraphs. \Cref{fig:graph} shows an instance where a geographic area, encoded by a graph containing $N=33$ nodes, is partitioned into $K=4$ territories, each of which is represented by a connected subgraph.
Thus, the spatial partitioning problem is equivalent to the graph partitioning problem (GPP) \cite{graphpartitioning} described next.

Given a positive integer $K \in {\mathbb{N}} \u{>1}$ and an undirected graph $\graph = \bigl( \vertexset, \edgeset \bigr)$ with non-negative edge weights, $\omega: \edgeset \rightarrow \mathbb{R}\u{>0}$, the solution to a GPP seeks a partition $\mathbf{\Pi}$ of $\graph$ with blocks of nodes $\mathbf{\Pi} = \bigl(V\u{1},\ldots,V\u{K}\bigr)$ such that $V\u{1} \cup \ldots \cup V\u{K} = \vertexset$ and $V\u{i} \cap V\u{j} = \phi \; \forall i,j \in \left\{1,2,\ldots,K\right\}, i \neq j$. 
 Alternatively, the output of a GPP can be represented by a plan on $\graph$ described by an assignment function $\xi : \vertexset \rightarrow \left\{1,2,\ldots,K\right\}$, where $\xi\bigl(v\bigr)=i$ implies that node $v$ is assigned to block $i$. A node $v$ is a neighbor of node $u$ if there is an edge $\bigl( u, v \bigr) \in \edgeset$. If a node $v \in V{\u i}$ has a neighbor $w \in V{\u j}$, $i \neq j$, then $v$ is called boundary node. Correspondingly, an edge that connects two boundary nodes is called \textit{cut edge} and $\edgeset \u{ij} = \left\{\left(u, v\right) \in \edgeset : u \in V\u{i}, v \in V\u{j}\right\}$ is the set of cut edges between two blocks, namely $V{\u i}$ and $V{\u j}$. 
An extra balance constraint may exist and enforces that all blocks have roughly equal weights, i.e., it requires that $\forall i \in \left\{1, \ldots, K\right\}: \left |V\u{i} \right | \leq \bigl( 1 + \tau \bigr) \left| \vertexset \right| / K$ for some threshold parameter $\tau \in \mathbb{R}\u{\geq0}$. Sometimes we also use weighted nodes with node weights. Weight functions on nodes and edges are extended to sets of such objects by summing their weights.
Note that a clustering is also a partition of the nodes. However, $K$ is usually not given in advance, and the balance constraint is removed. Note that a partition is also a clustering of a graph. In both cases, the goal is to minimize a particular objective called  the \textit{dispersion function}.
This is also identical to the connected $K-$partition problem, which partitions a graph into $K$ connected sub-graphs where $K \geq 3$, is a well-known \textsf{NP}-hard problem~\cite{dyer1985complexity}.

\subsection{Approaches}
\label{sec:meta}
In computer science, graphs constitute a preferred abstraction when modeling an application problem. 
Even if the application involves a different problem, partitioning a graph into smaller subgraph is often an important fundamental operation that helps to reduce problem complexity. Next, we briefly survey the different class of GP techniques broadly, an important end-application relevant to the spatial partitioning problems, and discuss the role of domain knowledge in designing spatial search techniques.
\begin{enumerate}
    \item \textit{Global algorithms} seeks a partition by directly working on the entire graph. Well-known global methods include \textit{exact algorithms}~\cite{BCR97,FMDS98,KRC00,Sen01,LR03,AFHM08,DW12} that rely on the branch-and-bound techniques~\cite{LD60}, \textit{spectral partitioning} techniques based on eigendecomposition of the Laplacian matrix, \textit{graph growing} approach based primarily on breadth-first search, \textit{flow-based methods} that make use of the max-flow min-cut theorem, and lastly the \textit{geometric partitioning} techniques that utilize the coordinates of graph nodes in space~\cite{Sim91,Wil91,GMT98}. The global algorithms are more suited to smaller graphs owing to high computation time, specially the exact methods. Also, these methods are mostly confined to bipartitioning but can be generalized to $k-$partitioning when $k$ is small.
    \item \textit{Iterative heuristics} start with an initial solution and tries to improve it through a variety of search operations. \textit{Local search} is the most widely used approach that updates a given solution by selecting a new one from the neighborhood. Different ways of defining the neighborhood and the selection strategy gives rise to a variety of techniques. Initial methods like the KL/FM method~\cite{KL70,FM82} was more suited for graph bisection. Later, $k-$way extensions to this method were proposed~\cite{KK98b,OS10}. Most existing local search algorithms  swaps nodes between adjacent blocks of the partition trying to minimize a dispersion function. This results in highly restrictive scope of possible improvement. For instance, the METIS approach cannot create balanced and contiguous partitions~\cite{ccnvd}.
    \item \textit{Multi-level approaches} perform graph partitioning by varying the granularity of the graphs~\cite{Wal04,Wal08,SS11}. It consists of the three main phases: coarsening, initial partitioning, and uncoarsening. Coarsening helps to reduce the problem size by iteratively approximating the original input graph with fewer degrees of freedom. This translates to substituting the parallel edges in the input graph with a single edge in the coarsened graph. Coarsening is terminated when the original graph is sufficiently small enough to be \textit{initially partitioned} using some (possibly expensive) exact methods discussed earlier. \textit{Uncoarsening} happens in two steps. Firstly, the partition in the coarse-level graph is translated back to a fine-level graph. Then, iterative improvement methods (discussed earlier) are usually applied to improve the partition. While multi-level approaches are successful for partitioning large graphs, it becomes challenging to tune this methods for graphs with fixed centers and highly varying degree of balance between the partitions.
    \item \textit{Metaheuristics} have been increasingly applied to the GP domain recently~\cite{KHKM11}.     
    There is a two-fold advantage of using metaheuristics. First, these frameworks are defined in a general sense, and hence can be modified to suit the needs of real-life optimization problems in terms of practical constraints like solution quality and execution time. Secondly, metaheuristics do not put any restriction on the optimization problem formulation (like constraints/ objective functions to be expressed as a linear function of the decision variables). Our focus is on \textit{population-based metaheuristics} like \textit{Evolutionary Algorithms} (EAs), which are derivative-free global optimization methods inspired by the process of natural evolution~\cite{eiben2003introduction}. An EA starts by initializing a population of trial solutions to the optimization problem, then it tries to improve the solutions via search operations, like recombination and selection, till a termination criteria is reached. In our work, we augment EAs with local search techniques for solving \spatial.
\end{enumerate}

\paragraph{Capacity-Constrained Network-Voronoi Diagram} Problems like districting, specially school districting, can be treated as a Capacity Constrained Network-Voronoi Diagram (CCNVD): ``\textit{Given a graph and a set of service center nodes, a CCNVD partitions the graph into a set of contiguous service areas that meet service center capacities and minimize the sum of the shortest distances from graph-nodes to allotted service centers}''~\cite{yang2013capacity,ccnvd}. For the school districting problem, the service centers represent the spatial units containing schools inside them. The Pressure Equalizer (PE) algorithm and its variants were proposed by~\citet{yang2013capacity,ccnvd} to address CCNVD. However, some important differences do exist. In the PE approach, the objective was to minimize the sum of the shortest distances from graph-nodes to their allotted service centers. Additionally, the following assumptions were made: all service centers have the same capacity, each non-service-center node has unit demand and all the service centers together could serve the demand of all the non-service-center nodes at any point of time. These simplistic assumptions may not apply for problems like districting. For instance, in school districting, the capacity of the schools and the student population corresponding to the graph nodes can vary considerably. Also, compactness is preferred to distance-based measures due to arbitrary shapes of spatial units forming a school district~\cite{biswas2020geospatial}.

\subsubsection{Hybrid metaheuristics}
Recently, researchers try to assimilate ideas from different classes of metaheuristics into a ``hybrid'' framework. One such framework, called \textit{memetic algorithm}~\cite{moscato1989evolution}, combines the local search technique with the recombination operator of EAs in order to balance exploration-exploitation. Thus, memetic algorithms benefit from the synergy between iterative improvement (exploitation) of the local search and the recombination operation (exploration) of the population-based methods. We take a step in this direction by integrating a randomized local search within a swarm-intelligent algorithm that mimics the foraging behavior of honeybees~\cite{akay2021survey}. 

However, adapting EAs to \spatial s is non-trivial due to the following:
\textit{Firstly}, EAs are designed to solve continuous-valued real-parameter global optimization problems. As such, they employ linear search moves for exploring the decision space by perturbing the incumbent solutions. This strategy is hardly suitable for discovering promising solutions in the discrete decision space of \spatial s.  
\textit{Secondly}, the presence of spatial constraints (topological properties), like contiguity, make \spatial s highly constrained in nature and harder to find feasible solutions.  
In fact, the infeasible solutions significantly outnumber the feasible solutions with an increase in problem size. This often results in expending tremendous computational effort in finding a feasible solution, especially when the search operators are not spatially cognizant.
\textit{Lastly}, EAs tend to reinitialize the solutions when they stagnate or violate constraints. In \spatial s, such a move is detrimental to preserving the goodness of solutions and lead to loss of valuable information. 

\subsubsection{Domain knowledge for spatial search}
In view of the above challenges, it is increasingly impracticable for vanilla EAs to solve \spatial s~\cite{hosage1986discrete}. This is mostly because the linear search moves are not suited to the discrete nature of problems we encounter. Besides, the traditional constraint handling techniques used in conjunction with EAs are of little help~\cite{mezura2011constraint}.
Hence, we use domain knowledge to guide the search process.
Domain knowledge refers to any auxiliary information that may enable a metaheuristic to efficiently search for feasible solutions. It includes both \textit{model-specific} and \textit{problem-specific} instructions.
Next, we discuss how model-specific domain knowledge is helpful in conducting spatially-aware search within an EA framework.

The first step in solving \spatial s is to instantiate initial feasible solution(s) and then improve the solution(s) locally by flipping spatial units between the adjacent (neighboring) territories~\cite{nagel1964simplified,openshaw1995algorithms}. There are two types of possible moves: a) \textit{move one unit} from its present (donor) territory to a neighboring (recipient) territory; b) \textit{swap/ exchange units} between two neighboring territories. The new solution produced by these moves are kept only if it is feasible and leads to an improvement in the objective function. Additionally, the local nature of the moves restricts the exploration of the decision space beyond the immediate neighborhood of the incumbent solution. However, these moves may also result in breaking the spatial contiguity of the territories involved in the move, thereby leading to an infeasible solution. Path relinking can be helpful in such scenarios for repairing the solutions if they enter the infeasible search space~\cite{glover2003scatter}. When infeasible solutions are made feasible again via repair operation, these solution(s) undergo strategic oscillations between the feasible and infeasible decision space, and may find better intermediate solutions~\cite{strategic}. 

 \section{Spatial partitioning as an optimization problem}
 \label{sec:optimization}
The optimization formulation corresponding to spatial partitioning problems can be written as

\begin{subequations}
\label{eq:pmedian}
    \begin{align}
        \mathbf{(P_0)}
        &
        \qquad 
        \overset{\mathbf{minimize}}{\mathbf{X} \in \mathcal{X}} \quad \objfunc \bigl( \mathbf{X} \bigr) = \sum\limits_{u= 1}^{N} \sum\limits_{v=1}^{N} X \u{uv} \cdot D\u{uv} \label{eq:pMPobjective}
        \\
        \text{s. t.}
        & 
        \nonumber
        \\
        &
        \sum\limits_{u = 1}^{N} X \u{uv} = 1,
        \qquad {} \qquad {}
        \qquad {} \qquad {}
        \quad {}
        \forall v = 1, 2, \ldots, N,
        \label{eq:pMPmutual}
        \\
        &
        \sum\limits_{u = 1}^{N} X \u{uu} = K,
        \qquad {} \qquad {}
        \qquad {} \qquad {}
        \qquad {}
        \label{eq:pMPtotal}
        \\
        &
        \left( 1 - \tau \right) {\mu} \cdot X \u {uu}  \leq \sum\limits_{v = 1}^{N} X \u{uv} \cdot {A \u v},
        \qquad {}
        \forall u = 1, 2, \ldots, N,
        \label{eq:pMPlowerbound}
        \\
        &
        \sum\limits_{v = 1}^{N} X \u{uv} \cdot {A \u v}  \leq \left( 1 + \tau \right) \capacity \cdot X \u {uu},
        \qquad {}
        \forall u = 1, 2, \ldots, N,
        \label{eq:pMPupperbound}
        \\
        &
        \sum\limits_{v\in {\bigcup}_{l \in S} \mathcal{N} ^l\backslash S} X \u{uv} - 
        \sum\limits_{v\in S} X\u{uv}  \ge 1-|S|,
        \forall u=1, 2,\ldots,N,
        \;  S \subset\{1, 2, \ldots, N\} 
         \backslash (\mathcal{N} ^u \cup \left\{ u \right\})
         \label{eq:pMPcontiguity}
          \\
        &
        X \u {uv} \in \bigl\{0, 1\bigr\},
        \qquad {}    \qquad {}    \qquad {}
        \qquad {}   \quad
        \forall u = 1, 2, \ldots, N,
        \;
        \forall v = 1, 2, \ldots, N.
        \label{eq:pMPbinary}
    \end{align}
\end{subequations}

where, 
$u$ and $v$ are the indices corresponding to the nodes of graph $\graph$;
$\mathbf{X} \in \left\{0, 1\right\} ^ {N \times N}$ is a binary assignment matrix, where $X\u{uu}=1$ implies that node $u$ is center node of a given subgraph;
$\mathbf{D} \in \mathbb{R} ^ {N \times N}$ is the distance (or dissimilarity) matrix, where $D\u{uv}$ is a distance\footnote{Normally, the Euclidean distance between the centroids of two spatial units is considered.} between nodes $u$ and $v$;
$\mathbf{A} \in \mathbf{R} ^{N \times 1} _{+}$ is the activity matrix, where $A\u v$ is an activity measure 
with respect to node $v$.
($\mathbf{P \u 0}$) is a constrained optimization problem with binary decision variables encoded by $\mathbf{\assignment}$.
Solving ($\mathbf{P} \u 0$) exactly will output an optimal partitioning plan encoded by the solution ${\mathbf{X}}^{*}$ that minimizes the objective function \eqref{eq:pMPobjective}, i.e., $\objfunc \bigl({\mathbf{X}}^{*}\bigr) \leq \objfunc \bigl( \mathbf{X} \bigr)$, $\forall \mathbf{X} \in \mathcal{X}$, subjected to a set of constraints \eqref{eq:pMPmutual}-\eqref{eq:pMPcontiguity}. $\mathcal{X}$ is the set of all possible partitioning plans, i.e., the assignment of $N$ spatial units to $K$ territories as show in \cref{fig:graph}. As spatial partitioning is analogous to graph partitioning, we shall interchangeably use the following groups of terms$-$spatial units/nodes and territories/subgraphs.

Constraint \eqref{eq:pMPmutual} enforces that each node is assigned uniquely to a subgraph. Constraint \eqref{eq:pMPtotal} ensures that the number of center nodes is $K$ implying that are exactly $K$ subgraphs since each subgraph has an unique center node. Constraints~\eqref{eq:pMPlowerbound}-\eqref{eq:pMPupperbound} restricts the total activity measure in a given subgraph to lie within a range of the mean activity measure $\mu$ (which is computed as $\mu = \frac{1}{K} \sum_{v=1}^{N} A \u v$) as measured by the tolerance parameter $\tau$. Usually the value of $\tau$ is kept in the range $\left[0.01, 0.1\right]$ depending on the application. These constraints are designed to ensure that a given subgraph has zero activity if it does not contain any center node, i.e., $X \u {uu}=0$. Lastly, the contiguity constraint~\eqref{eq:pMPcontiguity} that ensures that each subgraph is connected, where $\mathcal{N} ^u$ refers to all the nodes adjacent to a given node $u$. The connectivity constraint ensures that each territory is geographically contiguous, i.e., we can travel between any two points within a territory without crossing over to another adjacent territory. However, connectivity is expressed by an exponential number of comparisons, i.e., $\mathcal{O}\bigl( KN 2^N \bigr)$~\cite{drexl1999fast}, thereby making it impracticable to exactly solve moderate to large instances of this problem within a reasonable amount of time. 

\paragraph{Remarks} 
\Cref{eq:pmedian} can be solved by exact methods, like Linear Programming (LP), Integer Programming (IP), mixed-integer LP (MILP)~\cite{bertsimas1997introduction}. However, in trying to do so we made some interesting observations. \textit{Firstly}, these methods minimize a linear objective function based on dispersion, as in the $p-$median or the $p-$center problem~\cite{salazar2011new}. These linear measures of dispersion cannot account for optimizing the non-linear compactness metric. \textit{Secondly}, the exponentially big connectivity constraints may not guarantee territorial contiguity.
\textit{Lastly}, it may be difficult to find feasible solutions when the bound constraints \eqref{eq:pMPlowerbound} and \eqref{eq:pMPupperbound} are tightened by setting the value of $\tau$ to be sufficiently low.

\paragraph{Computational complexity}
Though the problem ($\mathbf{P \u 0}$) with exponential number of connectivity constraints, if we are given a $p-$partition $\graph' = \bigl( \vertexset', \edgeset' \bigr)$ of a graph $\graph = \bigl( \vertexset, \edgeset \bigr)$, we can check whether each subgraph of $\graph'$ is connected or not in polynomial time by using breadth-first-search algorithms. The feasibility of a solution to the problem can be verified in polynomial time, i.e., $\mathbf{P_0}$ is \textsf{NP}. Next, let us consider a particular instance where $\graph$ is a planar connected graph. If we take high values of tolerance parameter $\tau$, we can ensure that the balancing constraints \eqref{eq:pMPlowerbound} and \eqref{eq:pMPupperbound} are always satisfied. Then it becomes a $p$MP, which is a well-known \textsf{NP}-hard problem~\cite{hakimi1964optimum}. Since $p$MP is reducible to $\mathbf{P_0}$ in polynomial time, we can say that $\mathbf{P_0}$ is \textsf{NP}-hard. Interested readers may refer to \cite{dyer1985complexity} for an in-depth analysis of the computational complexity.


\subsection{School districting} 
\label{sec:school}
The school districting problem follows a spatial partitioning structure where the spatial units (or SPAs) are the nodes of a graph and the school boundaries (or SAZs) are the balanced, connected subgraphs. 
Hence, we reformulate \eqref{eq:pmedian} to define the school districting problem. We do make some adjustments in the optimization model based on the above-mentioned remarks and problem-specific domain knowledge. We shall visit them in turn.

Let a graph $\graph = \bigl(\vertexset, \edgeset \bigr)$ represent a school district with $N$ SPAs and $K$ schools. The number $K$ varies with the school level $\level:=$ ES, MS or HS, since we solve the problem at each level independently. 
Each node (SPA) can be represented as $A{\u v}=\bigl( {\u \mathrm{ES} }{p}{\u v}, {\u \mathrm{MS} }{p}{\u v}, {\u \mathrm{HS} }p{\u v}, {\u \mathrm{ES} }{c}{\u v}, {\u \mathrm{MS} }{c}{\u v}, {\u \mathrm{HS} }{c}{\u v} \bigr),\; v \in \left\{1, 2, \ldots, N \right\},$ where $\nX {\u v}$ is the student population residing in SPA $i$ corresponding to the school level~$\level$~(ES: elementary school, MS: middle school, HS: high school), and ${\u \level }{c}{\u v}$ is the program capacity of the schools contained in the same SPA. We assume that all the schools in a school district follow a consistent grade structure with respect to the school levels. 
For majority of the SPAs that don't enclose any school inside them, we have  ${\u \mathrm{ES} }{c}{\u v} = 0,$ ${\u \mathrm{MS} }{c}{\u v}=0$ and ${\u \mathrm{HS} }{c}{\u v}=0$. We consider a set of center nodes $\overline{\vertexset}=\left\{\overline{v}|X{\u {\overline{v} \overline{v}}} = 1 \right\}$, $\overline{\vertexset} \subset \vertexset$ corresponding to the SPAs containing schools inside them. Alternatively, we may write $\overline{\vertexset}=\{\overline{v}\u{1}, \overline{v}\u{2}, \ldots, \overline{v}\u{K}\}$, where $\overline{v}\u{i}$ is the node containing the $i$\textsuperscript{th} school. 

While drawing school boundaries, the following must be considered. \textit{Firstly}, each school has a different capacity to accommodate students. This is equivalent to the bound constraints \eqref{eq:pMPlowerbound} and \eqref{eq:pMPupperbound}, except that mean activity measure $\mu$ is replaced by the corresponding school's capacity ${\u \mathrm{L} }{c}{\u v}$. \textit{Secondly}, Euclidean distance between the centroid of the spatial units in \eqref{eq:pMPobjective} may not be a good representative of the commute time  due to the widely varying shapes of these units. Compactness measures that take into account the geometric shape can be a better alternative. These two considerations are linearly weighed using the weight factor $\lambda, 0 \leq \lambda \leq 1$, in formulating the optimization problem.
    \begin{align}
        \overset{\mathbf{minimize}}{\assignment \in \mathcal{X}}
            \mathcal{J}{\u s}\bigl( \assignment \bigr)
            =
            \lambda
            \underbrace{
            \sum\limits_{i=1}^{K}{
            \left | 1 - 
            \frac{
            \sum\limits_{u=1}^{N} \assign\u{u \overline{v}{\u i}} \cdot \u{\level}p{\u u}
            }{
            \sum\limits_{u=1}^{N} \assign\u{u \overline{v}{\u i}} \cdot \u{\level}c{\u u} 
            } 
            \right |
            }
            }_\text{target balance (aspatial)}
            \; + \;
            (1-\lambda)
            \underbrace{
            \sum\limits_{i=1}^{K}{
            \left |
            1 - 
            \frac{ 4 \pi \cdot
            \mathrm{Area} \bigl( \bigcup\u{u=1} ^{N} \{ u | {\assign}\u{ u \overline{v}{\u i} }=1\}\bigr)
            }{
            \left[\mathrm{Peri} \bigl( \bigcup\u{u=1} ^{N} \{ u | {\assign}\u{ u \overline{v}{\u i} }=1\}\bigr)\right]^2
            }
            \right |
                }
            }_\text{target compactness (spatial)}.
        \label{eq:objectivefunction}
    \end{align}
\newline
Several remarks are in order.
\textbf{1)} The \textit{target balance} measures the discrepancy between the schools' capacity and their attending student population. It attains the ideal value of $0$ when every school's attending student population (numerator) is equal to its program capacity (denominator). Note that target balance consolidates the bound constraints \eqref{eq:pMPlowerbound} and \eqref{eq:pMPupperbound} into an objective or a soft constraint. This helps to deal with solutions that cannot satisfy both these constraints in \eqref{eq:pmedian} by penalizing them heavily.
\textbf{2)} The \textit{target compactness} measures how far is a school's boundary from a perfectly compact shape (like a circle). We use the non-linear Polsby-Popper score~\cite{polsby1991third}, which is the ratio of the area of a zone to the area of a circle whose circumference is equal to the perimeter of the zone. The more compact the school boundaries become, the closer the value of target compactness gets to 0. 
\textbf{3)} Most of the school (re)districting happen to balance the student populations between existing schools. Hence, \textit{target balance} is given more weightage than \textit{target compactness}. In our setting, we always ensure that $\nicefrac{\lambda}{1 - \lambda} \geq 2$.
\textbf{4)} Prefixing the center nodes $\overline{\vertexset}$ helps to  satisfy the hard constraint \eqref{eq:pMPtotal} and automatically reduce the size of the optimization problem. Due to this advantage, we prefer to use problem-specific domain knowledge to perform prefixing. In the absence of any such information, clustering algorithms like $K-$medoids~\cite{kmedoids} can be used to determine a set of initial center nodes.
\textbf{5)} The remaining hard constraints, i.e., mutually exclusive assignment of nodes \eqref{eq:pMPmutual} and subgraph connectivity \eqref{eq:pMPcontiguity}, can be satisfied when initializing a feasible solution ${\assignment}$ and then perturbing it locally.
\textbf{6)} In solving \eqref{eq:objectivefunction}, the minimizer seeks a $K-$partition of a geographical area such that the territories are well-balanced, compactly-shaped and geographically contiguous. 
Overall, this approach generalizes to other spatial partitioning problems, like commercial territory design~\cite{rios2009reactive} and political redistricting~\cite{williams1995political}, that involve optimization of similar dispersion metric and subjected to constraints like balance, contiguity and compactness.

\section{The \algo~algorithm}
\label{sec:method}
\begin{figure}
    \centering
    \includegraphics[width=\linewidth, keepaspectratio]{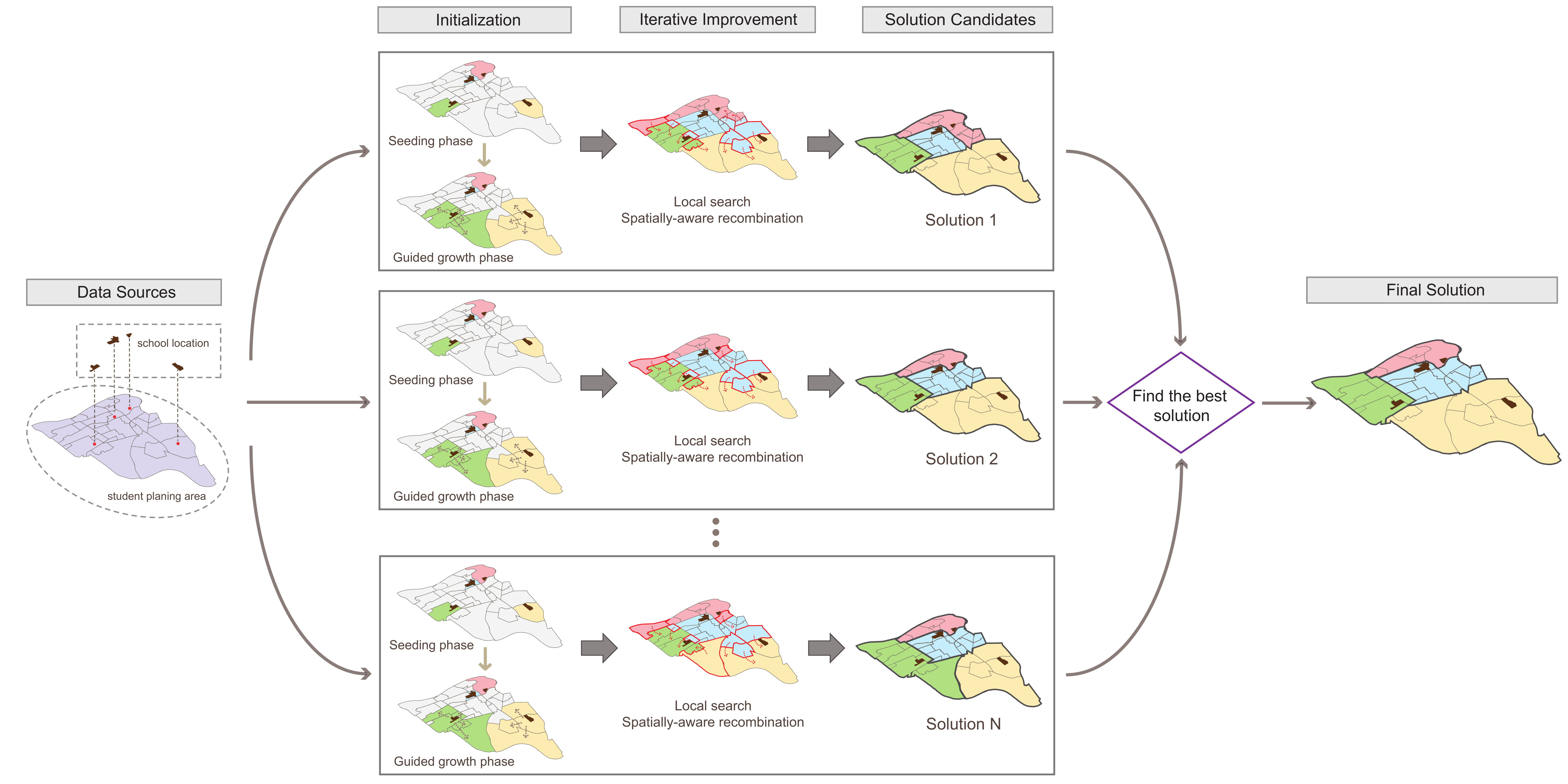}
    \caption{Outline of the \algo~approach for solving the school redistricting problem.}
\label{fig:algorithm}
\end{figure}

In this section, we describe our \algo~method for solving spatial partitioning problems. \algo~starts by initializing a population of randomly generated trial solutions as detailed in \cref{sec:initialize}. These solutions are iteratively improved till a termination criterion is met. The improvement takes place in two phases, \textit{local improvement} and \textit{spatially-aware recombination}, as detailed in~\cref{sec:local,sec:recombination}, respectively. The outline of our approach illustrated in \cref{fig:algorithm}.

\subsection{Initialization of trial solutions}
\label{sec:initialize}
\begin{figure}[htpt!]
        \centering
        \begin{subfigure}[b]{0.49\linewidth}
            \includegraphics[width=\textwidth, keepaspectratio]{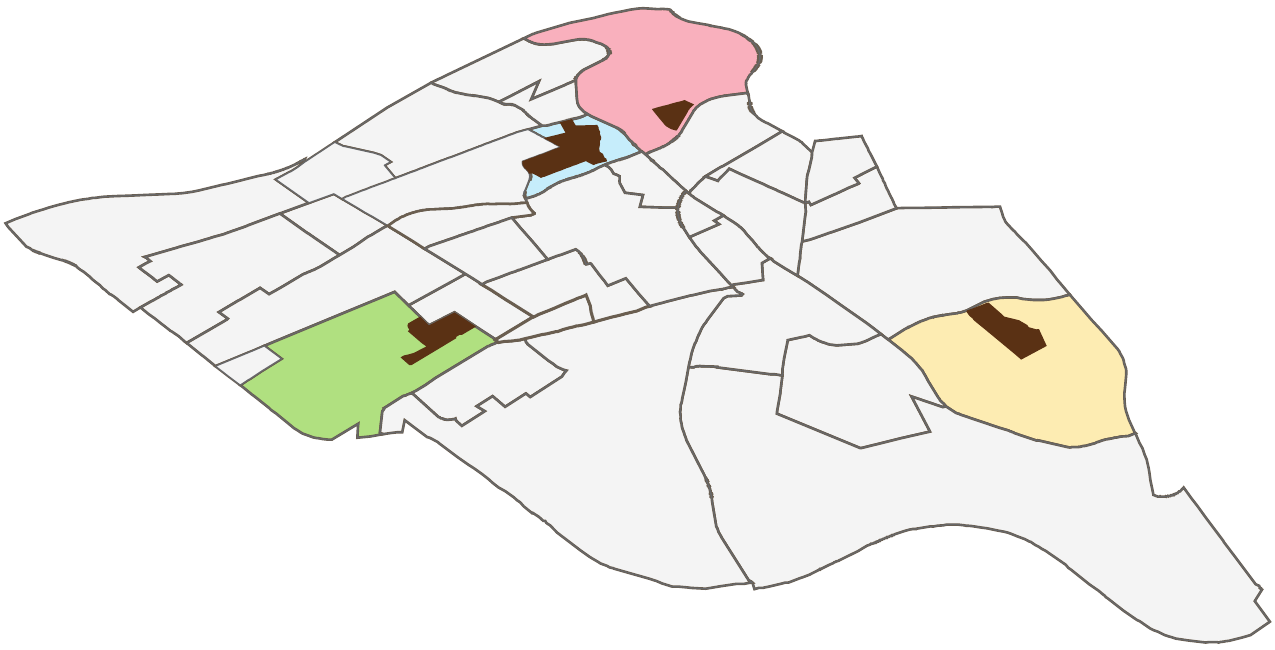}
            \caption{Seeding phase}
            \label{fig:seeding}
        \end{subfigure}
        ~\hfill~
       \begin{subfigure}[b]{0.49\linewidth}
            \includegraphics[width=\textwidth, keepaspectratio]{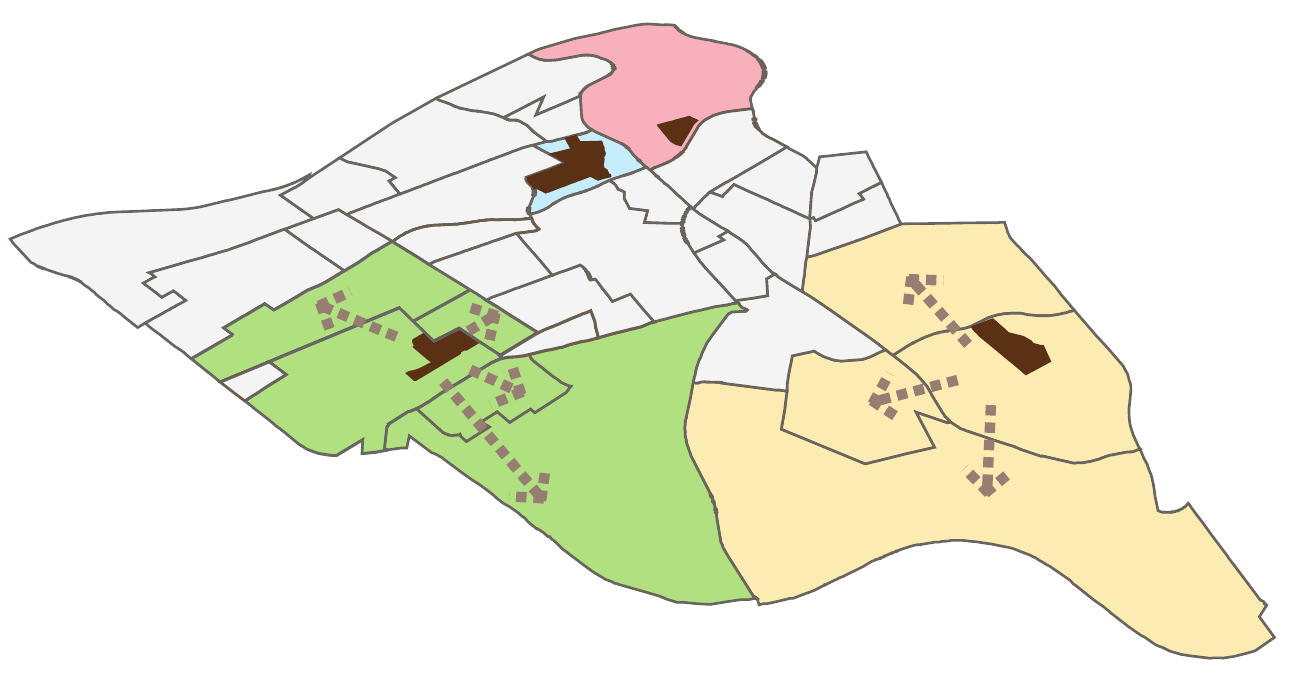}
            \caption{Guided growth phase}
            \label{fig:guided}
        \end{subfigure}
     \caption{The seeding phase (a) followed by the guided-growth phase (b) results in a new solution. Seeding identifies the spatial unit corresponding to the center nodes. The guided-growth phase helps to grow the territories by assigning the free spatial units (marked in light grey) based on the adjacency relation.}
    \label{fig:init-phase}
\end{figure}

Given a graph $\graph = \bigl( \vertexset, \edgeset \bigr)$, the initialization module instantiates a set of $Np$ trial solutions, $\mathcal{X} = \{{\assignment}^{(1)}, {\assignment}^{(2)}, \ldots, {\assignment}^{(Np)} \}$, where the $i$\textsuperscript{th} solution, ${\assignment}^{(i)}$, represents a particular partitioning of $\graph$. Note that we are overloading the notation on $\assignment$
The feasibility of these solutions are ensured by the  \textit{seeding phase} followed by the \textit{guided growth phase} as shown in \Cref{fig:init-phase}.
The pseudocode of initialization is provided in~\cref{algo:initialize}.

\paragraph{Seeding} This step helps to identify the seed units (prefix the center nodes) by leveraging problem-specific domain knowledge and assign each such unit (center node) to an unique territory (subgraph). In the context of school districting, seeding identifies each of the $K$ school-containing SPAs as center nodes as shown in~\cref{fig:seeding}. This leads to creation of $K$ partial subgraphs with just a single node in it, thereby ensuring that constraint \eqref{eq:pMPtotal} is satisfied. 
The assignment of these center nodes remain fixed throughout the partitioning process.

\paragraph{Guided growth} In the next step, the adjacency relationship between the spatial units are leveraged to grow the seed units into complete territories. This generates $K$ connected subgraphs representing a $K-$partition of $\graph$. \Cref{fig:guided} shows the guided growth phase, where a territory is randomly picked and grown by adding an adjacent (unassigned) spatial unit to it. The process is repeated till all the spatial units have been assigned to a territory, thereby satisfying constraints \eqref{eq:pMPmutual} and \eqref{eq:pMPcontiguity}.

\begin{algorithm}[htpt!]
\DontPrintSemicolon
\caption{Initialization}
\SetKwInOut{Input}{Input}
\SetKwInOut{Output}{Output}
\SetKwRepeat{Do}{do}{while}
\Input{Contiguity graph $\graph$, Population size $Np$, School level $\level$.}
 \Output{$\mathcal{X}:$ Population of trial solutions}
 \Begin{
 Determine the center nodes $\overline{\vertexset}$ for school level $\level$ and set $K \leftarrow | \overline{\vertexset} |$\;
 $\mathcal{X} \leftarrow \left\{ {} \right\}$ \Comment{Empty population}\\
\For{$i = \left\{1, 2, \ldots, Np \right\}$}{
$\vertexset$: Get the set of nodes in $\graph$\;
\textbf{Seeding phase} $\triangleright$\;
Set an initial assignment, i.e., ${\assignment}^{(i)} \in 0^{N \times N}$\;
\For{$\overline{v} \in \overline{\vertexset}$}{
    ${\assign} \u {\overline{v} \overline{v}} ^{(i)} \leftarrow 1$ \Comment{Assignment}\\
    $\vertexset \leftarrow \vertexset \backslash \left\{ \overline{v} \right\}$\;
}
\textbf{Guided-growth phase} $\triangleright$\;
\Do{$ \left| \vertexset \right| > 0$}{
     Randomly pick a center node $\overline{v},\, \overline{v} \in \overline{\vertexset}$\;
     Determine the subgraph $V$ containing $\overline{v}$, i.e., $V=\left\{u | u \in \vertexset,\,\assign \u {u \overline{v}} ^{(i)} = 1\right\}$ \;
     Find unassigned nodes adjacent to $V$, i.e., $\mathcal{N}\left( V \right)= \left\{ u | u \in \vertexset,\, \sum^{}_{v \in {\vertexset}}{{\assign} \u {u v} ^{(i)} = 0} \right\}$ \;
      \While{$ \left| \mathcal{N}\left( V \right) \right| > 0$ }
       {
        $u$: Randomly select a node from $\mathcal{N}\left( V \right)$\;
        $\assign \u {u \overline{v}} ^{(i)} \leftarrow 1$  \Comment{Assignment}\\
        $\mathcal{N}\left( V \right) \leftarrow \mathcal{N}\left( V \right) \backslash \left\{ v \right\}$ ,
        $\vertexset \leftarrow \vertexset \backslash \left\{ v \right\}$\;
            }
        }
     $\mathcal{X} \leftarrow \mathcal{X} \bigcup \left\{ {\assignment} ^{(i)} \right\}$\;
    }
}
\textbf{return} $\mathcal{X}$\;
\label{algo:initialize}
\end{algorithm}

 Note that during the growth phase, the adjacent spatial units are added in a random manner to grow the territories without consideration for the quality of the trial solutions generated. Hence, it is more than likely that the trial solutions will have low solution quality. To improve these solutions, we perform two steps of refinement discussed next.

\subsection{Local improvement}
\label{sec:local}

\begin{figure}
    \centering
    \includegraphics[keepaspectratio, width =0.3\textwidth]{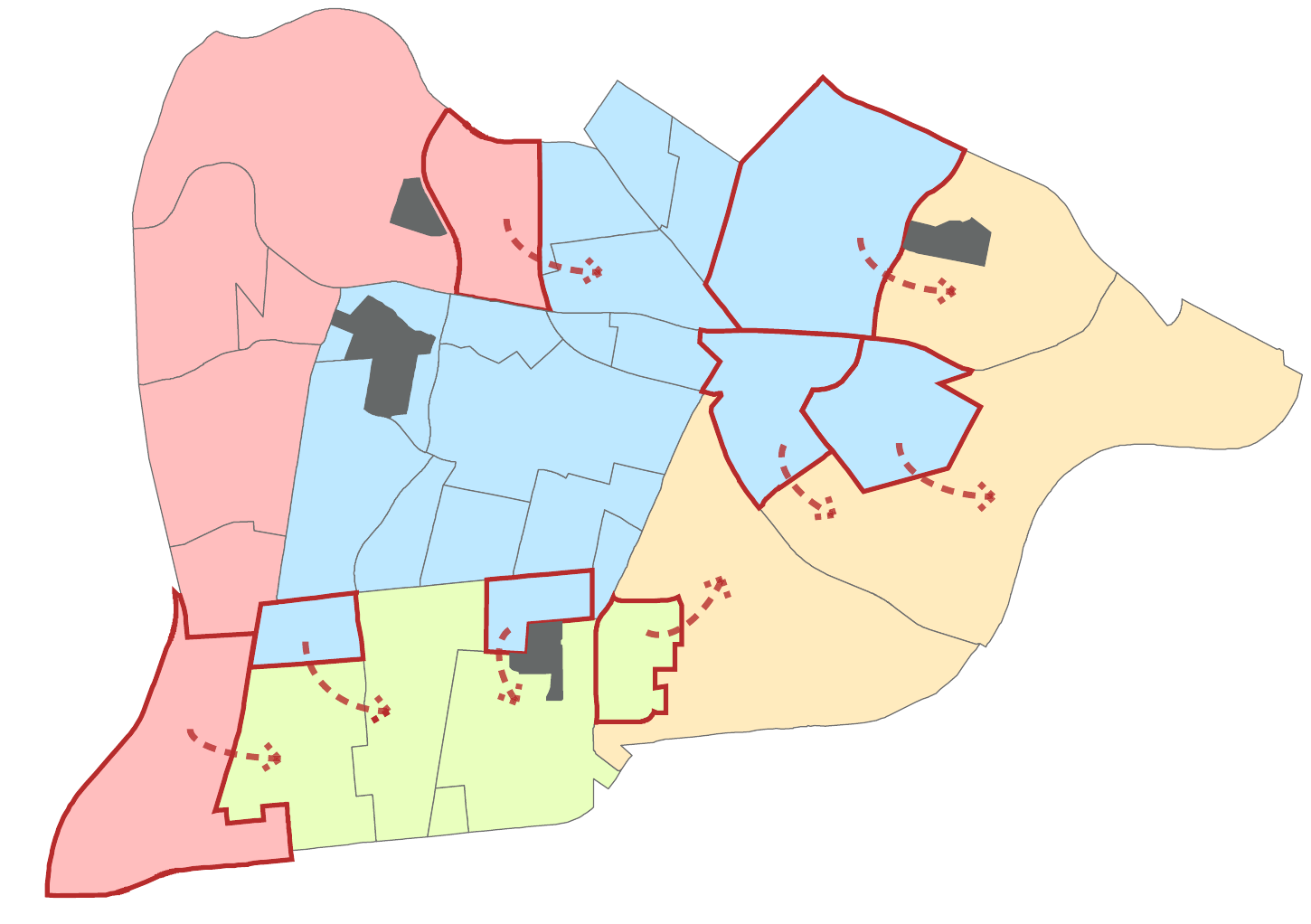}
    \includegraphics[keepaspectratio, width=0.65\textwidth]{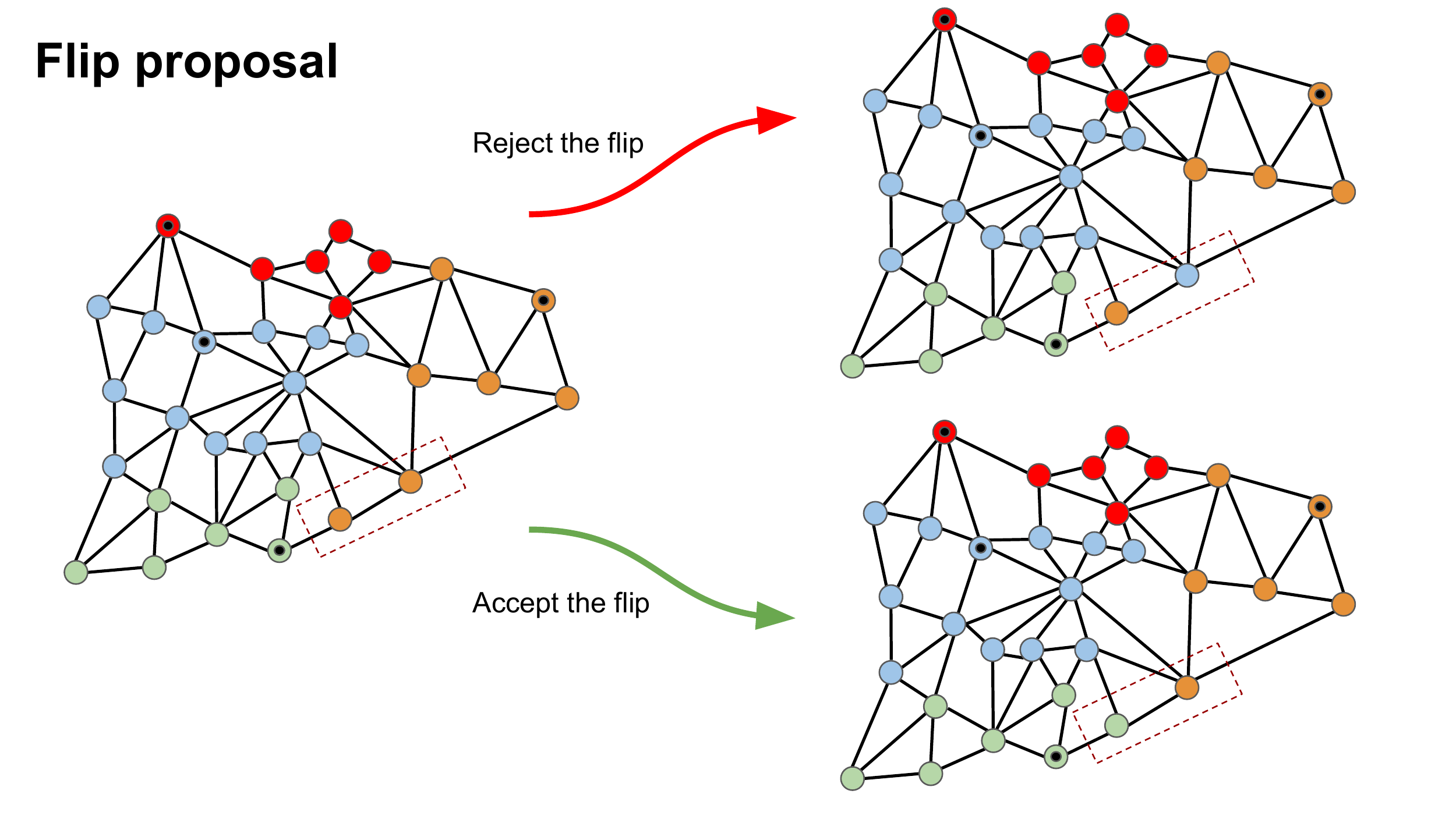}
        \caption{The local improvement helps to search for neighboring solutions that can be reached by flipping the membership of a boundary node. The \flip~proposal involves changing the assignment of a boundary node followed by acceptance/rejection of the new solution.}
        \label{fig:local}
\end{figure}

The local improvement searches the immediate neighborhood of an incumbent solution for improving the solution quality. 
Given the $i$\textsuperscript{th} solution, ${\assignment}^{(i)}$, we randomly pick a pair of subgraphs $V \u z$ and $V \u w$ ($w,z\in\left\{1, 2, \ldots, K\right\}$, $w \neq z$) such that they are adjacent, i.e., $|\edgeset \u{z,w}| > 0$. Then we may move a boundary node from $V \u z$ to $V \u w$ or vice-versa. This \textit{flipping} of nodes between adjacent subgraphs result in a new solution ${\Tilde{\assignment}}^{(i)}$. If this newly produced solution is of better quality, i.e., $\objfunc ( \Tilde{\assignment}{}^{(i)} ) < \objfunc ( {\assignment}{}^{(i)} )$, then $\Tilde{\assignment} ^{(i)}$ replaces ${\assignment} ^{(i)}$ in the population. Note that the connectivity of subgraphs $V\u z$ and/or $V \u w$ may be broken during flipping, thereby making $\Tilde{\assignment} ^{(i)}$ infeasible. 
To prevent such infeasibility, we only allow a move if it does not break the contiguity of the involved subgraphs. The local improvement operation is illustrated in~\cref{fig:local} and the pseudocode is provided in~\cref{algo:local}. Since the local improvement of a solution is independent of other solutions, it can leverage the parallel architecture of the computing platform. 

\begin{algorithm}
\DontPrintSemicolon
\caption{Local improvement}
\SetKwInOut{Input}{Input}
\SetKwInOut{Output}{Output}
\SetKwRepeat{Do}{do}{while}
\Input{Population of solutions $\mathcal{X}$, Contiguity graph $\graph$}
\Output{Updated solution}
\Begin{
\For{$i = \left\{1, 2, \ldots, Np\right\}$}
    {
    $\Tilde{\assignment} ^{(i)} \leftarrow {\assignment}^{(i)}$, $flipped \leftarrow False$\;
    \While{not flipped}{
        Randomly pick two adjacent subgraphs $V\u z$ and $V\u w$, i.e., $|\edgeset \u {z w}| > 0$\;
        Find boundary nodes in $V\u w$: $\mathcal{N}\left(V\u z\right) = \left\{v | (u, v)\in \edgeset\u {zw}, u \in V\u z \right\}$ \;
        \While{$|\mathcal{N}\left(V\u z\right)| > 0\;\&\&\; not\,moved$}{
        $v:$ From $\mathcal{N}\left( z \right)$ pick a random node $v$\;
        Move node $v$ from zone $V \u w$ to $ V \u z$, i.e.,
         $\Tilde{\assign}\u{zv} ^{(i)} \leftarrow 1$, $\Tilde{\assign}\u{zw} ^{(i)} \leftarrow 0$\;
        \If{$ V\u w \text{ and } V\u z \text{ are contiguous} $}{
                \eIf{${\objfunc (\Tilde{\assignment} ^{(i)} )} < \objfunc({\assignment}^{(i)}) \;\mathbf{||}\; rand\left(0,1\right) \leq p_r $}{
                ${\assignment}^{(i)} \leftarrow \Tilde{\assignment} ^{(i)},\,flipped\leftarrow True$   \Comment{Fitness-based replacement }\;
                }{
                $\Tilde{\assign}\u{zv} ^{(i)} \leftarrow 0$, $\Tilde{\assign}\u{zw} ^{(i)} \leftarrow 1$\Comment{Revert back the assignment}\;
                }
        }
        $\mathcal{N}\left(V\u z\right) \leftarrow \mathcal{N}\left(V\u z\right) \backslash \left\{v\right\}$
        }
    }

}
}
\label{algo:local}
\end{algorithm}
\normalsize

The random selection of subgraphs, i.e., $V\u z$ and $V\u w$, for performing node swaps may lead to redundancy. To prevent this, one may sequentially pick a subgraph, say $V\u z$, from a randomized list of subgraphs and determine adjacent subgraphs, say $V\u w$, for flipping the node. This is continued till a node flip is made. A flip is made when we find a better neighboring solution or we accept an inferior solution, i.e., ${\objfunc (\Tilde{\assignment} ^{(i)} )} > \objfunc({\assignment}^{(i)})$, with a very small probability $p_r$. While the former approach is greedy and prone to getting stuck at local optima, latter one  helps in randomization of the search move and is applied in metaheuristics like Simulated Annealing (SA)~\cite{kirkpatrick1983optimization}. 

\begin{figure}[htp]
    \begin{center}
        \includegraphics[keepaspectratio, width=0.6\textwidth]{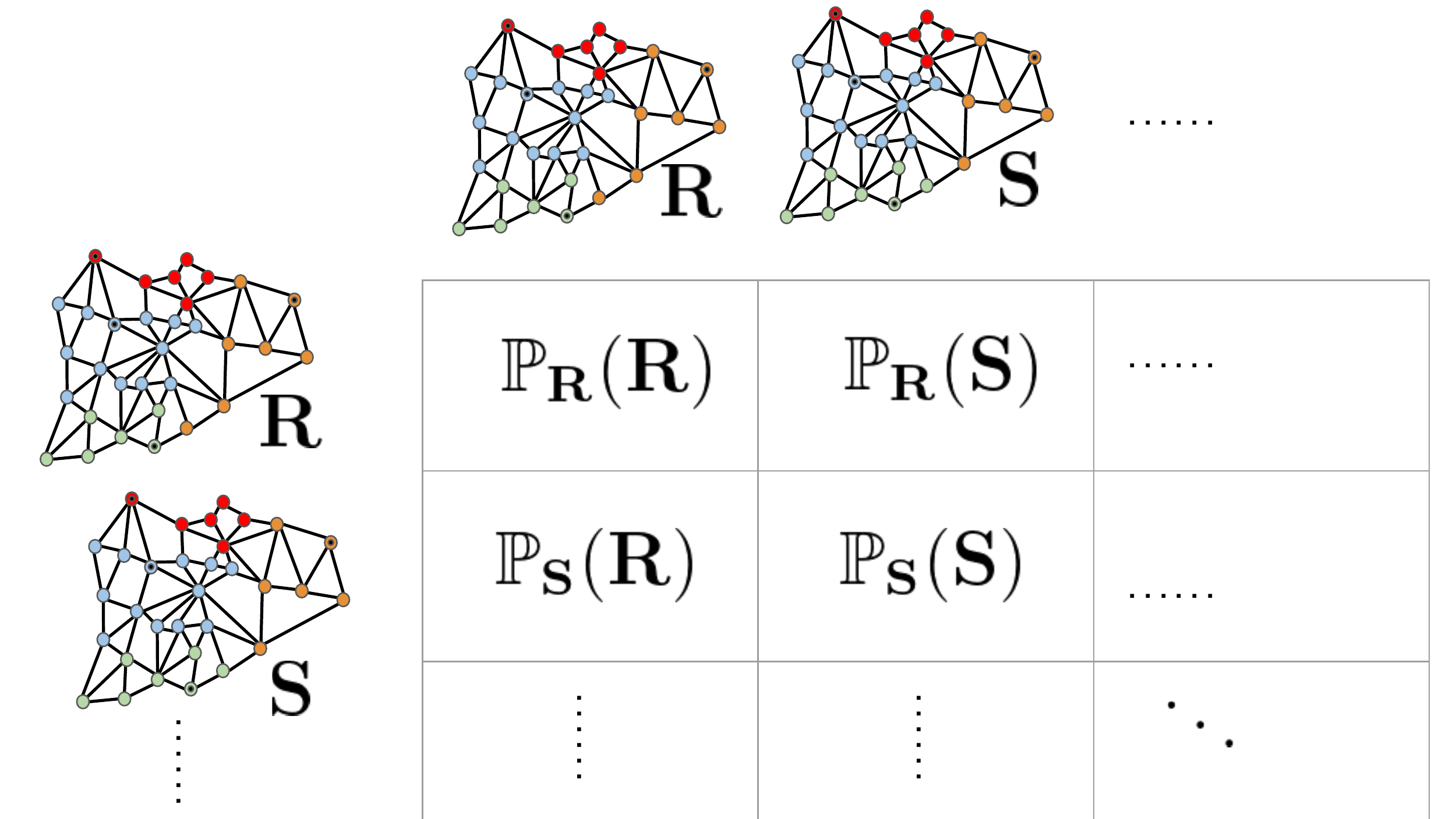}
    \caption{Theoretically if the flip proposal is carried out sufficient number of times, one may approximate the stationary distribution of transitions of the state space, where each state correspond to a districting plan.}
    \label{fig:states}
    \end{center}
\end{figure}

\paragraph{Markov Chains and Local Search}
 The local search mechanism here can be thought as instantiating the \flip-based walk, i.e., generating a new solution or \textit{districting plan} by changing the assignment of a single node as shown in \cref{fig:local}. Instantiating a series of flips to generate a sample of districting plans is akin to performing a random walk on the states of graph partitions and is encoded by a \flip-based Markov Chain~\cite{recombination}. Relatedly,  Markov chain Monte Carlo, popularly known as MCMC, is an effective technique for sampling owing to strong underlying theory, in the form of mixing theorems and convergence properties~\cite{mcmc}. In context of redistricting, lets imagine each districting plan representing a \textit{state} and a random walk is being performed on this state space. As the walker traverses from one state to another, we collect each state. On terminating the walk, this collection constitutes the representative sample of the plans. Performing the \flip-based walk involves changing the assignment of individual geographic units along district borders. In the standard MCMC paradigm, altering this basic step adjusts the stationary distribution. \cref{fig:states} gives a rough approximation of the idea.

The main purpose of these sampling-based techniques is to compare the a given districting plan in context of a \textit{representative sample}, i.e., a set of valid alternative plans. Closely following this is the need to relate the sampling distribution to the criteria set forth by domain experts. This may be a tough ask since any redistricting effort can be accompanied by a varying set of criteria, some of which are difficulty to quantify objectively. Our objective here is different. We use a customized sampling distribution to generate an ensemble of plans and save the best quality plan, as determined by an objective function.

\subsection{Spatially-aware recombination}
\label{sec:recombination}
During local improvement, the individual solutions are improved independently without any exchange of information between them. 
Interestingly, it is possible to determine a better intermediate solution by combining features from two solutions. Population-based methods enable mixing of solutions through the recombination operation~\cite{eiben2003introduction}. This results in better exploration of the search space. However, the vanilla recombination operation is not suitable for spatially-constrained problems. Hence, we develop the \textit{spatially-aware recombination} operation.

The recombination operation is motivated by the exchange of genetic material between different organisms which leads to production of offspring. In this process, two (parent) solutions, say ${\assignment}^{(i)}$ and ${\assignment}^{(j)}$, are selected such that ${\assignment}^{(i)}$ is picked randomly while solution ${\assignment}^{(j)}$ is selected probabilistically based on the fitness value. The fitness function is defined to allow solutions with lower functional value have higher fitness as this is a minimization problem. For maximization, the fitness can be set equal to the objective functional value. We expect that ${\assignment}^{(j)}$ is fitter than ${\assignment}^{(i)}$ and thus  ${\assignment}^{(i)}$ can learn from  ${\assignment}^{(j)}$. The steps of recombination operation are provided in~\cref{algo:recombination}. 

\begin{algorithm}
    \DontPrintSemicolon
    \caption{Spatially-aware recombination}
    \SetKwInOut{Input}{Input}
    \SetKwInOut{Output}{Output}
    \SetKwRepeat{Do}{do}{while}
    \Input{$\mathcal{X}:$ Population of solutions, $\graph:$ Contiguity graph}
    \Output{Updated solution}
    \Begin{
    Find the fitness values:
    $\mathbf{\mathcal{H}} ^{(i)} =  \left.\frac{1}{1 + | \objfunc \left({\assignment}^{(i)} \right) |} \right| \; \forall i = 1, 2, \ldots, |\mathcal{X}|$\;
    \For{$i = \left\{1, 2, \ldots, |\mathcal{X}|\right\}$}
        {$\Tilde{\assignment} ^{(i)} \leftarrow {\assignment}^{(i)}$, $\Tilde{\assignment} ^{(j)}$: Probabilistically selected $j^{th}$ solution based on the fitness value\;
        Randomly pick a subgraph $V$ such that $ 0 < \left|V^{(i)} \cap V^{(j)} \right| < \min \left(|V^{(i)}|, |V^{(j)}|\right)$\;
        Find the set of incoming nodes
        $I{\u V} = \left\{v | v\in V^{(j)} \setminus V^{(i)} \text{ and } \exists u \in V^{(i)} \text{ s.t. } (u,v) \in \edgeset \right\}$ and outgoing nodes
        $O{\u V} = \left\{u | u\in V^{(i)} \setminus V^{(j)} \text{ and } \exists v \in V^{(j)} \text{ s.t. } (u,v) \in \edgeset \right\}$\;
        Randomly pick an incoming node $v \in I \u V$ and an outgoing node $u \in O \u V$ \;
        Simultaneously insert node $v$ into $V^{(i)}$ and remove node $u$ from zone $V^{(i)}$; also update the assignments in $\Tilde{\assignment} ^{(i)}$\;
        If $V^{(i)}$ has rendered non-contiguous by the swap operation, repair $\Tilde{\assignment} ^{(i)}$\;
    }
    \For{$i = \{1, 2, \ldots, |\mathcal{X}|\}$}
        {
        \If{${\mathcal{F} \left(\Tilde{\assignment} ^{(i)}\right)} \leq {\mathcal{F}} \left({\assignment}^{(i)}\right)$}{
            ${\assignment}^{(i)} \leftarrow \Tilde{\assignment} ^{(i)}$ \Comment{Fitness-based update}
            }
        }
        
    }
    \label{algo:recombination}
    \end{algorithm}
\normalsize

\begin{figure}
    \centering
    \includegraphics[keepaspectratio, width =0.95\linewidth]{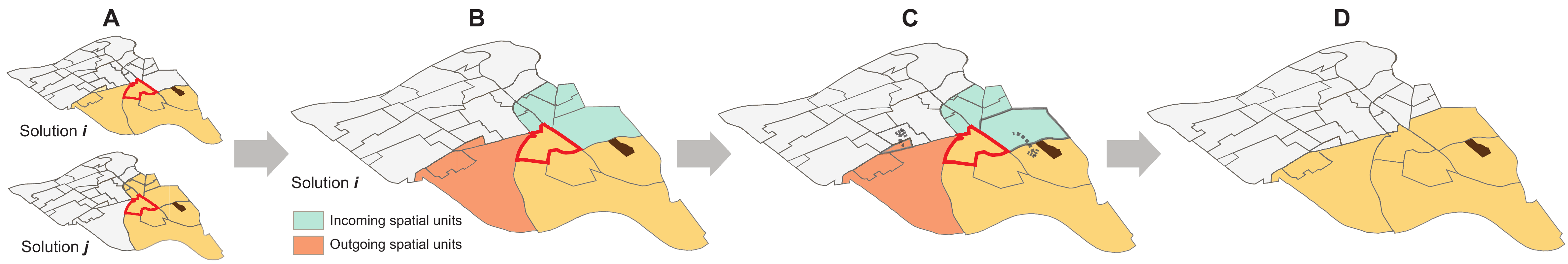}
    \caption{Illustrating the individual steps involved in the spatially-aware recombination operator.}
    \label{fig:recombination}
\end{figure}

Suppose a subgraph $V$ is present in both solutions $i$ and $j$, marked as $V^{(i)}$ and $V^{(j)}$, such that they have a common node. Every subgraph should satisfy this condition since the center nodes remain unchanged. The subgraph $V^{(i)}$ is modified by simultaneously inserting a node $v$ (present in $V^{(j)}$ but not in $V^{(i)}$) and deleting a node $u$  (present in $V^{(i)}$ but not in $V^{(j)}$). This swapping of node steers solution ${\assignment}^{(i)}$ towards the fitter solution ${\assignment}^{(j)}$ as illustrated in Figure~\ref{fig:recombination}. In doing so, we expect to find intermediate solutions that may have better fitness than the incumbent solution ${\assignment}^{(i)}$.

\begin{figure}
    \centering
    \includegraphics[width=0.95\linewidth, keepaspectratio]{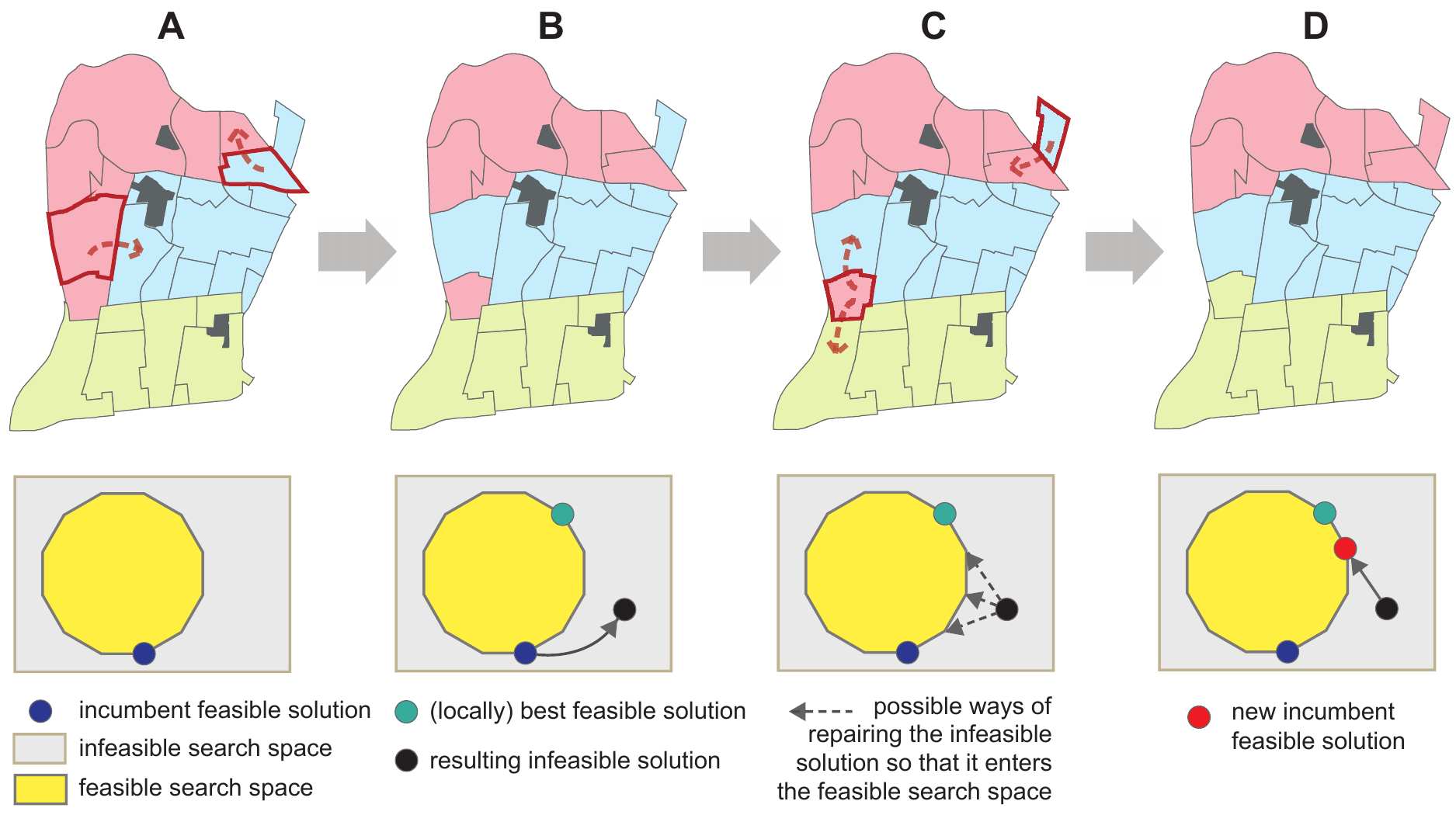}
    \caption{During the recombination, the swapping of spatial units (A) may result in an infeasible solution (B), which needs to be repaired (C) to make the solution feasible (D). As a result of these moves, the solution oscillates between feasible and infeasible search space as demonstrated above. For illustration purposes, we have depicted this oscillation through a continuous search space. We notice that the \textsf{recombination} plus \textsf{repair} operation may result in a solution that differs from the incumbent solution by multiple hops. As a  result, it may approach (locally) optimal solutions.}
    \label{fig:repair}
\end{figure}

Interestingly, the swapping of the nodes may break the connectivity of the involved subgraphs. To reduce the chances of such undesirable scenarios, we perform the swap operation using boundary nodes. Nevertheless, a repair operation still needs to be applied in case the connectivity of the subgraphs are broken.
To repair a solution, we use the breadth-first search (BFS) traversal for enumerating the connected components in the disconnected subgraph, say $V$. Then, each connected component is analyzed for the presence of the center node. If center node is absent, all the nodes in this component is reassigned to the neighboring subgraphs. When no prior information about the center nodes is available, we may retain the largest-sized connected component of $V$ and reassign the other components. The repaired solution $\Tilde{\assignment} ^{(i)}$ might be few steps away from the incumbent solution ${\assignment} ^{(i)}$ in discrete space and thus helps in controlled exploration of the search space. The advantage of repair operation is shown in \cref{fig:repair}. Note that we have used a depiction of continuous search space in \cref{fig:repair} though this is a discrete optimization problem. This was done to simply show the movement of a solution through the repair process.

 The solutions newly generated from the recombination operation needs to be updated in the population synchronously based on their fitness values. The solution update is important to keep the fitter solutions in the population so as to make the search progress. With careful implementation, this step can be parallelized in modern computer systems by using multiprocessing.

\section{Experimentation}
\label{sec:exp}
In this section, we conduct experiments with the proposed \algo~method and some well-known metaheuristics on the school districting problem.  The dataset description is provided in~\cref{sec:dataset} followed by the performance metrics in~\cref{sec:evaluation}. Details about setting the objective function and the baseline methods are provided in \cref{sec:objectivefunction,sec:baselines} respectively. The comparative evaluation of the baseline methods is performed in~\cref{sec:comparison} followed by studies on the effect of solution initialization in \cref{sec:initialization} and computational analysis in 
\cref{sec:computational}. Next, \cref{sec:ablation} focuses on a series of ablation tests for investigating the individual components of \algo. Lastly, \cref{sec:planning} studies how \algo~can aid planners in real-life planning.

\subsection{Dataset}
\label{sec:dataset}
The study was performed on two school districts (counties) located in the mid-Atlantic region of the USA. These school districts have seen recent population growth in certain areas, thereby making the problem challenging for the \spatial s tested here. The following GIS data attributes of both the district swere used for experimentation.
\begin{itemize}
    \item \textbf{SPAs}: The location coordinates of the spatial units along with aggregated student count at different school levels (Elementary, Middle and High).
    \item \textbf{Schools}: The location of the school building, school level, and its program capacity.
\end{itemize}

In comparison to the dataset used in \cite{spatial}, the only difference is that here we used the data for the  school year $2020-2021$. We resorted to using the new dataset since many parents unenrolled their children from the public schools at the onset of the COVID-19 pandemic resulting in more population imbalance. This would present more challenging problem scenarios for the redistricting algorithms. Table~\ref{ch5:tab:data} presents the summary statistics of this new dataset.
  
\begin{table}
\centering
  \caption{Summary statistics of the school districts for the school year 2020-21.}
  \label{ch5:tab:data}
  \begin{tabular}{c|c|ccc}
    \toprule
    \multirow{2}{*}{District} & \multirow{2}{*}{\#SPAs ($N$)}  & \multicolumn{3}{c}{\#Schools ($K$)}
    \\ \cline{3-5}
    & & Elementary  & Middle & High
    \\ \hline
    \lcps &  453 &  57 & 16 & 16
    \\ 
    \fcps & 1313 & 138 & 26 & 24\\ 
    \bottomrule
  \end{tabular}
\end{table}

 We performed a few additional pre-processing steps. Specifically, we modeled the SPAs as nodes in a graph and generated the adjacency relationship between the nodes. Also, we determined the center nodes, i.e., the nodes corresponding to the spatial units containing schools inside them, by performing point-in-polygon test using the PySal library~\cite{rey2010pysal}.

\subsection{Evaluation metrics}
\label{sec:evaluation}
The solution to the school redistricting problem generated by an algorithm is actually a plan or a zoning configuration of the school boundaries. For evaluating the plan, we utilized two performance metrics that can be interpreted as percentage scores since they lie in the range $\left[0, 100\right]$.

\begin{itemize}
    \item \textbf{Balance} measures the average balance between a school's program capacity and the number of students residing within its boundary. It is calculated as
    \begin{equation}
    \text{bal}\left(\assignment\right) = 
    100 \times
    \left| 1 -
     \frac{1}{K}
     \sum\limits_{i=1}^{K}{
            \left | 1 - 
            \frac{
            \sum\limits_{u=1}^{N} \assign\u{u \overline{v}{\u i}} \cdot \u{\level}p{\u u}
            }{
            \sum\limits_{u=1}^{N} \assign\u{u \overline{v}{\u i}} \cdot \u{\level}c{\u u} 
            } 
            \right |
            }
    \right|
    \end{equation} 
    We penalized both under-enrolled and overburdened schools equally with respect to the capacity of schools. This is an important metric for school planners since most of the boundary changes occur to achieve a better balance in schools.
    
    \item \textbf{Compactness} measures how tightly a school's boundary is packed on an average with respect to its perimeter. A scaled Polsby-Popper score~\cite{polsby1991third} is used to measure compactness as
    \begin{equation}
    \text{com}\left(\assignment\right) = 
     \frac{100}{K}
     \sum\limits_{i=1}^{K}{
            \left |
            \frac{ 4 \pi \cdot
            \mathrm{Area} \bigl( \bigcup\u{u=1} ^{N} \{ u | {\assign}\u{ u \overline{v}{\u i} }=1\}\bigr)
            }{
            \left[\mathrm{Peri} \bigl( \bigcup\u{u=1} ^{N} \{ u | {\assign}\u{ u \overline{v}{\u i} }=1\}\bigr)\right]^2
            }
            \right |
                }
    \end{equation} 
    Compact school boundaries often translate to proximal schools that students can walk to and thereby lower the transportation cost incurred by the school district.

\end{itemize}

\subsection{Setting the objective function}
\label{sec:objectivefunction}
The two objectives- balance and compactness - are conflicting in nature~\cite{biswas2019regal}. Hence an optimizer is needed to trade-off between the objectives while solving the problem.  We ensured that the condition $\nicefrac{\lambda}{1 - \lambda} \geq 2$ held for both the districts while defining the objective function in~\cref{eq:objectivefunction}, i.e., we set the weighing parameter $\lambda$ in \cref{eq:objectivefunction} to $0.7$ and $0.8$ for districts~\lcps~and~\fcps~, respectively. This could ensure that the importance of population balance is at least twice\footnote{Balancing a school's attending population w.r.t. its program capacity is the driving forces behind school redistricting. This is consistent with the actual practices of school planners, who give more importance to balance over compactness.} as that of compactness. The value of $\lambda$ was higher for district $\mathbf{Y}$ since it had higher population imbalance in the schools to begin with. 
For each school district, we independently solved three instances of the redistricting problem$-$ Elementary School (ES), Middle School (MS), and High School (HS). We observed varying characteristics of the problem in each instance.

The ES instance of the problem is more challenging than the others due to different factors. Firstly, a school district that has seen recent population growth is most likely due to an influx of young children, which often leads to burgeoning demand for new ESs. Oftentimes, the new ESs are situated at arbitrary locations without sufficient separation between them. This goes against the clustering assumption of well-separated cluster centers. Secondly, the ESs exhibit a wide variation in their program capacity, with the newly built ones having higher capacity than their older counterparts, thereby making it difficult to balance the student population with the schools' capacity. In attempting to satisfy the schools' capacities, the optimization algorithm may fill in concave segments in the school boundaries with regular-shaped spatial units having a high density of student population. The MS and HS, being well-separated and showing fewer deviations of capacity, are comparatively easier to solve. Interestingly, the ES boundaries are more compact than their MS and HS counterparts. In comparison to ES SAZs, a greater proportion of boundaries of MS and HS share borders with the school district's boundary, which is usually zigzagged by naturally occurring geographies (highly irregular geometries). These additional considerations, besides the spatial constraints, make school redistricting a challenging \spatial~to solve.

\subsection{Baseline methods}
\label{sec:baselines}

The following baseline methods are used for comparative study:

\begin{itemize}
    \item \textbf{Local search-based techniques}: We set the parametric configuration for each of these baselines based on the literature~\cite{biswas2019regal}.

    \begin{itemize}
    
        \item \textsf{Stochastic Hill Climbing (SHC)}~\cite{SHC}: A variation of the basic Hill Climbing that searches the immediate neighborhood of a feasible solution in a random manner. If an equally good or a better solution is found, it replaces and the search continues till a local optimum is obtained.
        
        \item \textsf{Simulated Annealing (SA)}~\cite{kirkpatrick1983optimization}: A stochastic version of the Hill Climbing that is based on the process of tempering of metals. It allows for worsening moves to take place if no better solutions are found and can escape local optima.
        
        \item \textsf{Tabu Search (TS)}~\cite{glover1998tabu}:  An algorithm uses a restrictive (tabu) list to forbid revisiting recently explored solutions so that the new neighboring solutions can be explored.
        
    \end{itemize}
    
    \item \textbf{Sampling-based techniques}: We include three sampling-based optimization techniques developed based on the link between MCMC and redistricting stated earlier in \Cref{sec:local}. 
    \begin{itemize}
        \item \textsf{B}alanced, \textsf{A}lways \textsf{A}ccept (\modelA)
        \item \textsf{B}alanced and \textsf{C}ompact, \textsf{A}lways \textsf{A}ccept (\modelB)
        \item \textsf{A}ccept \textsf{I}mproving \textsf{O}bjective (\modelC)
    \end{itemize}
    For more details about these techniques, refer to Chapter 5 in \cite{biswas2022phd}.
    
    \item \textbf{SPATIAL}~\cite{spatial}: The population size was set to 10 and 20 for districts \lcps~and~\fcps, respectively. Trial runs were simulated till 1000 and 2000 iterations for districts \lcps~and~\fcps , respectively.
    
\end{itemize}

Note that we considered two more baseline methods$-$Greedy Randomized Adaptive Search Procedure (GRASP)~\cite{rios2009reactive} and Mixed Integer Linear Programming (MILP)~\cite{tshirabe}. While GRASP's performance was inferior to the other baselines, MILP could not converge to a feasible solution for 4/6 test cases even with a run-time budget of 24 hours. For codes of SA, TS, SHC and SPATIAL, go to \url{https://github.com/subhodipbiswas/SpatialPartitioning} and for the sampling-based techniques, check out \url{https://github.com/subhodipbiswas/SamplingbasedSchoolRedistricting}.

\subsection{Comparison with existing methods}
\label{sec:comparison}

For comparison purposes, we simulated $25$ trial runs of each baseline and recorded the final solutions. Each solution represented a districting plan of school boundaries, which were evaluated based on the metrics defined in Section~\ref{sec:evaluation}. In Table~\ref{tab:comparison}, we reported the mean and standard deviation of these metrics.
We also included the existing school boundary configuration of the school districts. It is marked as \texttt{Existing} in the table. The results revealed that~\algo~was able to generate better quality solutions in the majority of the test cases. This especially held for district~\lcps. Besides achieving better balance, \algo~obtained improved compactness score. The difference is specially marked when you compare the \texttt{Existing} plans with the plans generated by \algo.

\begin{table}
    \centering
    \caption{Performance of peer algorithms on the problem of school boundary formation in both the districts. The best performing entries are marked in boldface.}
    \label{tab:comparison}
    \setlength{\tabcolsep}{1pt}
    \footnotesize
    \begin{tabular}{@{}c|cc|cc|cc@{}}
    \toprule
    \multicolumn{7}{c}{District \lcps}
    \\ 
    \hline
    \multirow{2}{*}{
    Models
    }  & 
    \multicolumn{2}{c}{Elementary School} &
    \multicolumn{2}{c}{Middle School} &
    \multicolumn{2}{c}{High School}
    \\ \cline{2-7} &
        Balance     &    Compactness    &
        Balance     &    Compactness    &
        Balance     &    Compactness    \\ 
        \midrule
        \texttt{Existing} &
        83.5020 $\pm$ 0.0000 & 32.5344 $\pm$ 0.0000 &  
        89.7379 $\pm$ 0.0000 & 26.7671 $\pm$ 0.0000 & 
        87.0786 $\pm$ 0.0000 & 27.3452 $\pm$ 0.0000   
    \\
    \midrule
    SA & 
    87.7697$\pm$0.8280 & 38.1032$\pm$1.6977 & 
    92.3789$\pm$0.5726 & 32.3574$\pm$3.6421 & 
    96.4240$\pm$1.9741 & 26.9094$\pm$2.7749 
    \\
    TS & 
    87.3788$\pm$0.7079 & 36.4537$\pm$1.5931 & 
    92.5729$\pm$0.2888 & 33.0756$\pm$2.1707 & 
    95.9435$\pm$1.8527 & 28.3494$\pm$2.1524 
    \\
    SHC &
    86.6755$\pm$0.8642 & 36.5780$\pm$1.7189 & 
    92.5583$\pm$0.1806 & 32.3115$\pm$2.3647 & 
    95.6461$\pm$1.9090 & 27.9708$\pm$2.1571 
    \\
    \midrule
    \modelA &
    71.7328$\pm$2.3193 & 30.5329$\pm$1.5741 & 
    89.4501$\pm$1.1493 & 18.6148$\pm$1.6085 & 
    90.3080$\pm$1.8968 & 16.8370$\pm$1.7453 
    \\
    \modelB &
    71.0820$\pm$2.4890 & 30.7455$\pm$1.3047 & 
    90.1128$\pm$1.5444 & 18.6181$\pm$1.7215 & 
    90.9247$\pm$1.6280 & 17.0758$\pm$2.0001 
    \\
    \modelC & 
    86.6930$\pm$1.0247 & 37.4850$\pm$1.5981 & 
    92.4021$\pm$0.4920 & 33.4105$\pm$2.3518 & 
    95.7360$\pm$1.7613 & 30.5854$\pm$2.3411 
    \\
    \midrule
    \algo &
    \textbf{87.9353$\pm$0.6175} & \textbf{38.8988$\pm$1.4759} & 
    \textbf{92.5926$\pm$0.1191} & \textbf{37.5398$\pm$1.6283} & 
    \textbf{97.7948$\pm$0.4602} & \textbf{31.5193$\pm$1.8568}  
    \\
    \bottomrule
    \multicolumn{7}{c}{District \fcps} 
    \\ 
    \hline
    \multirow{2}{*}{
    Models
    }  & 
    \multicolumn{2}{c}{Elementary School} &
    \multicolumn{2}{c}{Middle School} &
    \multicolumn{2}{c}{High School}
    \\ \cline{2-7} &
        Balance     &    Compactness    &
        Balance     &    Compactness    &
        Balance     &    Compactness    \\ 
        \midrule
    \texttt{Existing} &
        82.3835 $\pm$ 0.0000 & \textbf{35.9234 $\pm$ 0.0000} &  
        84.2310 $\pm$ 0.0000 & 27.7096 $\pm$ 0.0000 & 
        86.9541 $\pm$ 0.0000 & 26.8006 $\pm$ 0.0000   
    \\
    \midrule
    SA & 
    94.6386$\pm$1.0634 & 30.0782$\pm$1.4172 &
    91.3987$\pm$1.1006 & 22.9419$\pm$2.4456 & 
    93.1795$\pm$1.5845 & 24.2182$\pm$2.8013
    \\
    TS & 
    93.3600$\pm$0.7183 & 29.7212$\pm$0.6997 &
    92.1146$\pm$0.4123 & 25.5295$\pm$2.8529 &
    \textbf{93.6894$\pm$1.4057} & 26.4599$\pm$2.7025 
    \\
    SHC & 
    92.9656$\pm$0.9595 & 29.7697$\pm$0.9929 &
    91.4258$\pm$0.8386 & 24.5064$\pm$2.5644 &
    93.2757$\pm$1.4755 & 24.3142$\pm$2.1289 
    \\
    \midrule
    \modelA &
    67.4731$\pm$1.9623 & 27.1551$\pm$0.8379 & 
    86.0361$\pm$1.2005 & 10.3140$\pm$1.1395 &
    84.8948$\pm$2.1824 & 10.1874$\pm$0.8616
    \\
    \modelB &
    67.3870$\pm$2.1147 & 27.2899$\pm$0.7668 &
    86.6818$\pm$1.4137 & 10.6439$\pm$1.1396 &
    85.4392$\pm$1.7517 & 10.1476$\pm$0.7313 
    \\
    \modelC & 
    93.0592$\pm$0.8565 & 30.7064$\pm$0.7959 &
    91.9071$\pm$0.5007 & 25.7025$\pm$2.6570 &
    93.4927$\pm$1.4288 & 26.1186$\pm$2.2065 
    \\
    \midrule
    \algo &
    \textbf{94.9097$\pm$0.4351} & 30.2618$\pm$0.9105 & 
    \textbf{92.3780$\pm$0.1295} & 27.8174$\pm$2.2644 & 
    92.6337$\pm$0.6078 & \textbf{29.6600$\pm$1.5477}
    \\
      \bottomrule
    \end{tabular}
    \normalsize
\end{table}

For district~\fcps, the performance of \algo~was comparable to the baseline in terms of the balance scores, yet the compactness scores were comparatively better. \algo~was the only model that achieves at par or better compactness than the \texttt{Existing} plan. We observed similar trend for district~\lcps. The baseline methods adopted a greedy approach by continuing to look for better solutions in the local neighborhood of the incumbent solution. In doing so, the solutions lying just outside their immediate neighborhood remained elusive to them. On the other hand, the spatially-aware recombination technique enabled~\algo~to find intermediate solutions outside the immediate neighborhood. The repair operation was particularly instrumental in finding such solutions, some of which may be better than the solutions presented in the immediate neighborhood. This is in line with the findings reported in~\cite{strategic}. 

The sampling-based methods, with the exception of \modelC, do not result in high quality plans since they are not optimizing on a particular objective. They are designed to operate like random search methods without a greedy selection procedure. \modelC, on the other hand, applies greedy selection procedure and in essence similar to the techniques like SA, TS and SHC.

\subsection{
Effect of solution initialization
}
\label{sec:initialization}
To inspect if solution initialization has any effect on the quality of solutions, we created a variant of the baseline methods indicated by an asterisk$-$\textsf{Algorithm}$^*$ is the newly created variant of \textsf{Algorithm}. If say \algo~starts with randomly generated initial solutions (as depicted in \cref{sec:initialize}),  the initial solutions in \algo$^*$ will correspond to the \textsf{existing} boundary configuration of the school district under consideration. Since the school districts redraw their boundaries frequently (in response to changing needs), the \textsf{existing} boundary configuration represents a (locally) good solution to the problem. This is akin to solving the \textit{school redistricting} problem where we redraw the school boundaries instead of designing them from scratch as in school districting. This subtlety is more of a matter of technicality. Interestingly, while performing school redistricting, we ensured that each school boundary was geographically contiguous to being with. If not, we would use the repair operation outlined earlier to reassign some SPAs for creating contiguous school boundaries.

We ran 25 simulations on each school level of both the districts and tabulated the performance metrics in \cref{tab:initialization}. We observe that on an average the \textsf{Algorithm}$^*$ variants are better than \textsf{Algorithm} in both the metrics. This is mostly because the \texttt{existing} solution that \textsf{Algorithm}$^*$ starts with is of better quality than the randomly-generated solutions used by \textsf{Algorithm}. For district \lcps, there is a clear trend showing that \textsf{Algorithm}$^*$ is better than \textsf{Algorithm} across all the possible metrics and problem instances. However, we do see some exceptions in district \fcps, especially for MS and HS problem instances. This can be attributed to the \textsf{Algorithm}$^*$ variants getting stuck in a local optima. This is highly plausible since the initial solutions of  \textsf{Algorithm}$^*$ are very similar\footnote{
The only difference between the solutions are due to the reassignment of discontinuous SPAs of school boundaries for the purposes of maintaining the contiguity of the school boundaries.
} to each other. On the other hand, \textsf{Algorithm} may manage to escape the local optima courtesy the widely varying range of the starting solutions it is initialized with. The number of local optima increases exponentially with the increase in problem size. Hence, this trend is expected since district \fcps~is almost three times the size of district \lcps. Note that both the variants of \algo~were able to achieve superior results in majority of the cases.

\begin{table}
    \centering
    \caption{The effect of solution initialization on peer algorithms.  The best performing entries are marked in boldface.}
    \label{tab:initialization}
    \setlength{\tabcolsep}{1pt}
    \footnotesize
    \begin{tabular}{@{}c|cc|cc|cc@{}}
    \toprule
    \multicolumn{7}{c}{District \lcps}
    \\ 
    \hline
    \multirow{2}{*}{
    Models
    }  & 
    \multicolumn{2}{c}{Elementary School} &
    \multicolumn{2}{c}{Middle School} &
    \multicolumn{2}{c}{High School}
    \\ \cline{2-7} &
        Balance     &    Compactness    &
        Balance     &    Compactness    &
        Balance     &    Compactness    \\ 
        \midrule
    \texttt{Existing} &
        83.5020 $\pm$ 0.0000 & 32.5344 $\pm$ 0.0000 &  
        89.7379 $\pm$ 0.0000 & 26.7671 $\pm$ 0.0000 & 
        87.0786 $\pm$ 0.0000 & 27.3452 $\pm$ 0.0000   
    \\
    \midrule
    SA & 
    87.7697$\pm$0.8280 & 38.1032$\pm$1.6977 & 
    92.3789$\pm$0.5726 & 32.3574$\pm$3.6421 & 
    96.4240$\pm$1.9741 & 26.9094$\pm$2.7749 
    \\
    SA$^*$ & 
    87.9809$\pm$0.5044 & 40.9826$\pm$1.0907 &
    92.5345$\pm$0.2965 & 37.1108$\pm$2.1572 &
    97.3246$\pm$0.5107 & 33.0759$\pm$1.8799 
    \\
    \hline
    TS & 
    87.3788$\pm$0.7079 & 36.4537$\pm$1.5931 & 
    92.5729$\pm$0.2888 & 33.0756$\pm$2.1707 & 
    95.9435$\pm$1.8527 & 28.3494$\pm$2.1524 
    \\
    TS$^*$ & 
    88.2290$\pm$0.3177 & 40.4481$\pm$0.5532 &
    92.7145$\pm$0.0019 & 38.1273$\pm$0.1339 & 
    97.3869$\pm$0.1536 & 32.7663$\pm$0.3416 
    \\
    \hline
    SHC &
    86.6755$\pm$0.8642 & 36.5780$\pm$1.7189 & 
    92.5583$\pm$0.1806 & 32.3115$\pm$2.3647 & 
    95.6461$\pm$1.9090 & 27.9708$\pm$2.1571 
    \\
    SHC$^*$ & 
    88.0565$\pm$0.3398 & 40.4219$\pm$0.6413 &
    92.6503$\pm$0.1009 & 37.9347$\pm$1.1750 &
    97.5288$\pm$0.6344 & 33.3857$\pm$1.6935
    \\
    \midrule
    \modelA &
    71.7328$\pm$2.3193 & 30.5329$\pm$1.5741 & 
    89.4501$\pm$1.1493 & 18.6148$\pm$1.6085 & 
    90.3080$\pm$1.8968 & 16.8370$\pm$1.7453 
    \\
     \modelA$^*$  & 
     83.0738$\pm$0.6002 & 32.4948$\pm$1.0189 &
     92.1247$\pm$0.3935 & 28.4179$\pm$0.8139 &
     94.5031$\pm$0.7808 & 26.3284$\pm$1.3810 
    \\
    \hline
    \modelB &
    71.0820$\pm$2.4890 & 30.7455$\pm$1.3047 & 
    90.1128$\pm$1.5444 & 18.6181$\pm$1.7215 & 
    90.9247$\pm$1.6280 & 17.0758$\pm$2.0001 
    \\
    \modelB $^*$ &
    84.1403$\pm$0.5138 & 32.2777$\pm$1.0458 & 
    92.2637$\pm$0.3897 & 29.6189$\pm$1.0030 & 
    95.5958$\pm$0.8448 & 26.2145$\pm$1.6023 
    \\
    \hline
    \modelC & 
    86.6930$\pm$1.0247 & 37.4850$\pm$1.5981 & 
    92.4021$\pm$0.4920 & 33.4105$\pm$2.3518 & 
    95.7360$\pm$1.7613 & 30.5854$\pm$2.3411 
    \\
    \modelC $^*$ &
    87.7953$\pm$0.4931 & 41.0278$\pm$0.7990 &
    92.6865$\pm$0.0594 & 37.1659$\pm$0.8954 &
    96.6412$\pm$1.0957 & 33.7013$\pm$1.7627 
    \\
    \midrule
    \algo &
    87.9353$\pm$0.6175 & 38.8988$\pm$1.4759 & 
    92.5926$\pm$0.1191 & 37.5398$\pm$1.6283 & 
    97.7948$\pm$0.4602 & 31.5193$\pm$1.8568  
    \\
    \algo$^*$ &
    \textbf{88.2474$\pm$0.3008} & \textbf{41.6610$\pm$0.9948} &
    \textbf{92.6639$\pm$0.0642} & \textbf{39.9138$\pm$0.7408} &
    \textbf{98.0286$\pm$0.2052} & \textbf{35.2453$\pm$1.2159}
    \\
    \bottomrule
    \multicolumn{7}{c}{District \fcps} 
    \\ 
    \hline
    \multirow{2}{*}{
    Models
    }  & 
    \multicolumn{2}{c}{Elementary School} &
    \multicolumn{2}{c}{Middle School} &
    \multicolumn{2}{c}{High School}
    \\ \cline{2-7} &
        Balance     &    Compactness    &
        Balance     &    Compactness    &
        Balance     &    Compactness    \\ 
        \midrule
    \texttt{Existing} &
        82.3835 $\pm$ 0.0000 & 35.9234 $\pm$ 0.0000 &  
        84.2310 $\pm$ 0.0000 & 27.7096 $\pm$ 0.0000 & 
        86.9541 $\pm$ 0.0000 & 26.8006 $\pm$ 0.0000   
    \\
    \midrule
    SA & 
    94.6386$\pm$1.0634 & 30.0782$\pm$1.4172 &
    91.3987$\pm$1.1006 & 22.9419$\pm$2.4456 & 
    93.1795$\pm$1.5845 & 24.2182$\pm$2.8013
    \\
    SA\textsuperscript{*} &
    95.1757$\pm$0.5959 & 33.9204$\pm$0.9149 &
    91.4492$\pm$0.6952 & 28.2519$\pm$1.6041 &
    92.2339$\pm$0.4015 & 30.4304$\pm$1.2035 
    \\
    \hline
    TS & 
    93.3600$\pm$0.7183 & 29.7212$\pm$0.6997 &
    92.1146$\pm$0.4123 & 25.5295$\pm$2.8529 &
    \textbf{93.6894$\pm$1.4057} & 26.4599$\pm$2.7025 
    \\
    TS\textsuperscript{*} & 
    93.9882$\pm$0.2143 & \textbf{37.4500$\pm$0.4358} &
    89.8580$\pm$0.6495 & \textbf{31.9823$\pm$1.0187} &
    92.0073$\pm$0.1227 & 34.3255$\pm$0.7516
    \\
    \hline
    SHC & 
    92.9656$\pm$0.9595 & 29.7697$\pm$0.9929 &
    91.4258$\pm$0.8386 & 24.5064$\pm$2.5644 &
    93.2757$\pm$1.4755 & 24.3142$\pm$2.1289 
    \\
    SHC\textsuperscript{*} & 
    94.2043$\pm$0.5056 & 35.3233$\pm$0.6058 &
    90.7112$\pm$1.3864 & 29.9091$\pm$1.3810 &
    91.9250$\pm$0.0950 & 32.1575$\pm$1.0527 
    \\
    \midrule
    \modelA &
    67.4731$\pm$1.9623 & 27.1551$\pm$0.8379 & 
    86.0361$\pm$1.2005 & 10.3140$\pm$1.1395 &
    84.8948$\pm$2.1824 & 10.1874$\pm$0.8616
    \\
    \modelA\textsuperscript{*} &
    80.1400$\pm$0.4253 & 32.3159$\pm$0.4638 &
    81.5796$\pm$1.0063 & 22.8711$\pm$0.9850 &
    87.7556$\pm$0.6411 & 22.6101$\pm$0.7538 
    \\
    \hline
    \modelB &
    67.3870$\pm$2.1147 & 27.2899$\pm$0.7668 &
    86.6818$\pm$1.4137 & 10.6439$\pm$1.1396 &
    85.4392$\pm$1.7517 & 10.1476$\pm$0.7313 
    \\
    \modelB\textsuperscript{*} &
    81.3120$\pm$0.3927 & 32.0952$\pm$0.5303 &
    83.1391$\pm$0.9270 & 22.4642$\pm$0.7627 &
    88.9047$\pm$0.5496 & 22.2048$\pm$0.7405 
    \\
    \hline
    \modelC & 
    93.0592$\pm$0.8565 & 30.7064$\pm$0.7959 &
    91.9071$\pm$0.5007 & 25.7025$\pm$2.6570 &
    93.4927$\pm$1.4288 & 26.1186$\pm$2.2065 
    \\
    \modelC\textsuperscript{*} &
        94.3157$\pm$0.5046 & 35.8568$\pm$0.7109 &
        91.3974$\pm$0.9010 & 31.1680$\pm$1.5590 &
        91.8972$\pm$0.1067 & 33.0105$\pm$1.1274 
    \\
    \midrule
    \algo &
    94.6537$\pm$0.3691 & 30.1714$\pm$0.6381 &
    \textbf{92.3885$\pm$0.1194} & 28.9142$\pm$1.1682 &
    92.6819$\pm$0.7202 & 29.5414$\pm$1.3540
    \\
    \algo\textsuperscript{*} &
    \textbf{95.3731$\pm$0.3617} & 36.1492$\pm$0.7468 &
    91.7854$\pm$0.2029 & 31.8754$\pm$1.3194 &
    91.9185$\pm$0.1156 & \textbf{34.6683$\pm$0.9838}
    \\
      \bottomrule
    \end{tabular}
    \normalsize
\end{table}

Nevertheless, the advantage of \textsf{Algorithm}* variants like \algo* is in reconfiguration of the existing school boundaries, i.e., it is useful in scenarios like opening of a new school or closure of an existing school. However, to enable \textsf{Algorithm}* work successfully, it is desired that the existing plan has a high percentage of connected subgraphs. For instance, some of the existing plans of district~\fcps~did not have geographic contiguity (subgraph connectivity) and hence the disconnected subgraphs have to be repaired. We noticed that a higher proportion of disconnected subgraphs, on being repaired, results in arbitrarily-shaped districts. For such districts, it may be difficult to arrive at a better configuration due to the local structure that the problem imposes.

\subsection{Computational complexity and ideas on scaling up}
\label{sec:computational}
\begin{figure}
    \centering
    \begin{subfigure}[b]{\linewidth}
        \includegraphics[width=0.495\textwidth, keepaspectratio]{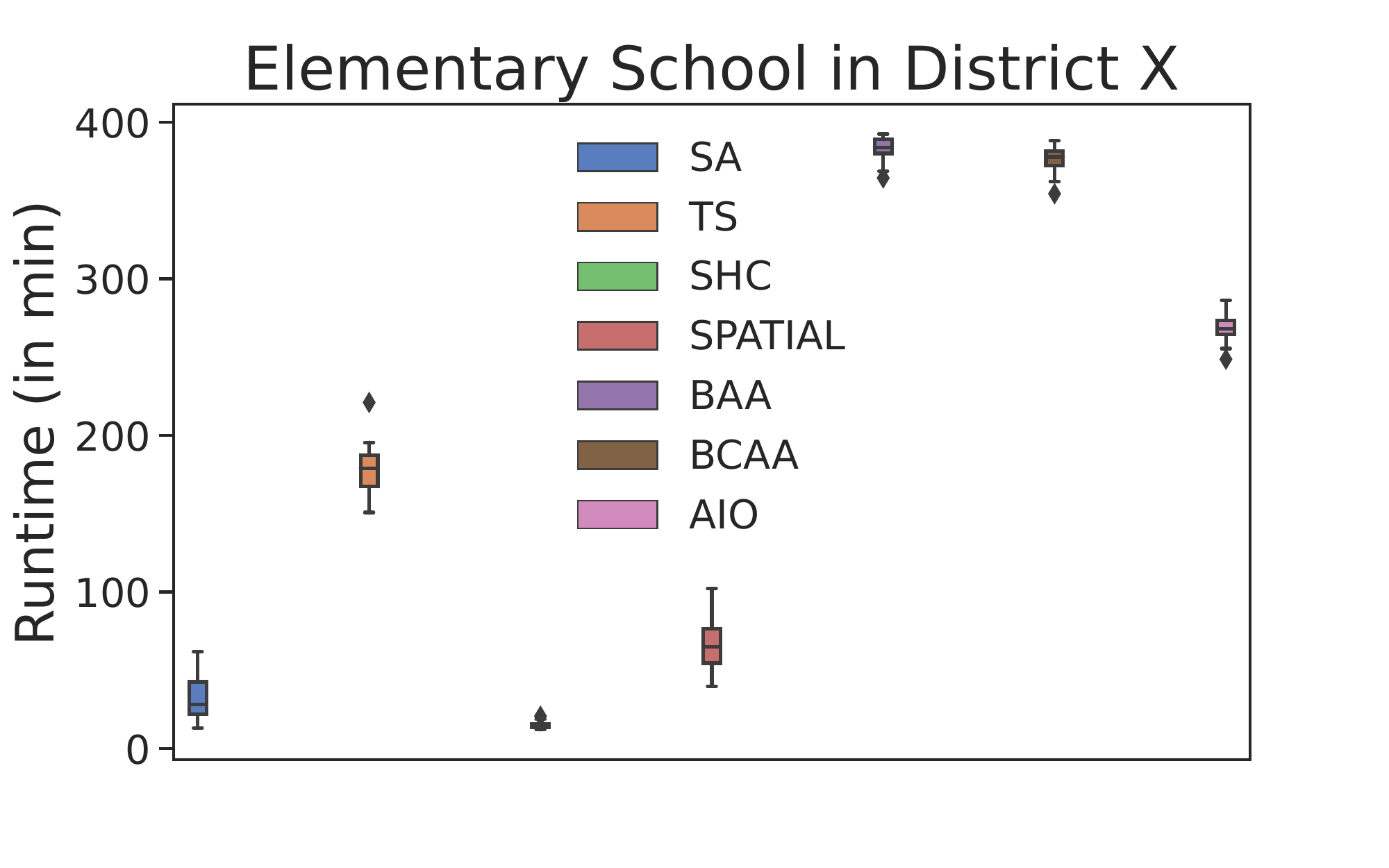}
    ~    
        \includegraphics[width=0.495\textwidth, keepaspectratio]{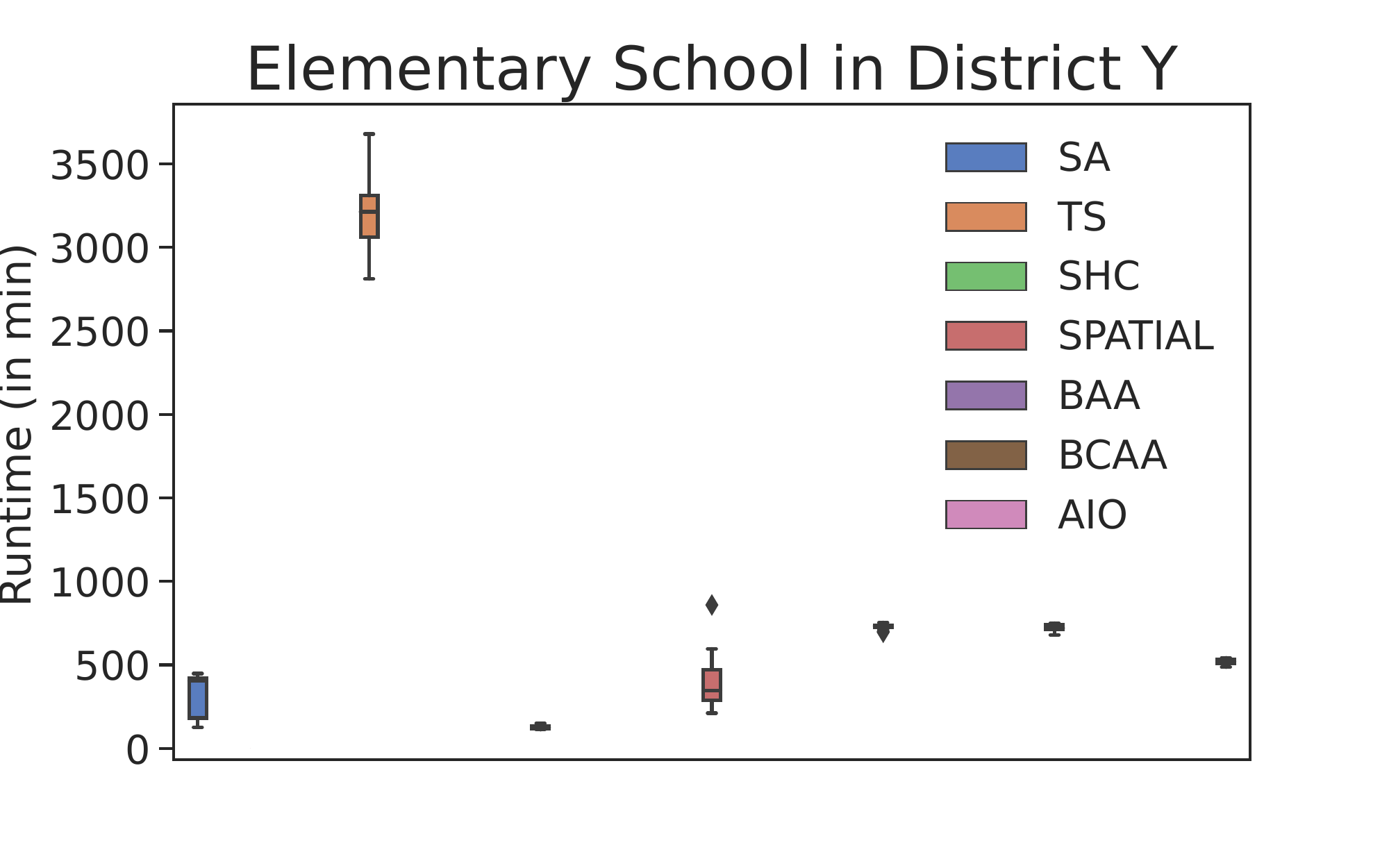}
        \caption{Elementary school}
    \end{subfigure}

   \begin{subfigure}[b]{\linewidth}
        \includegraphics[width=0.495\textwidth, keepaspectratio]{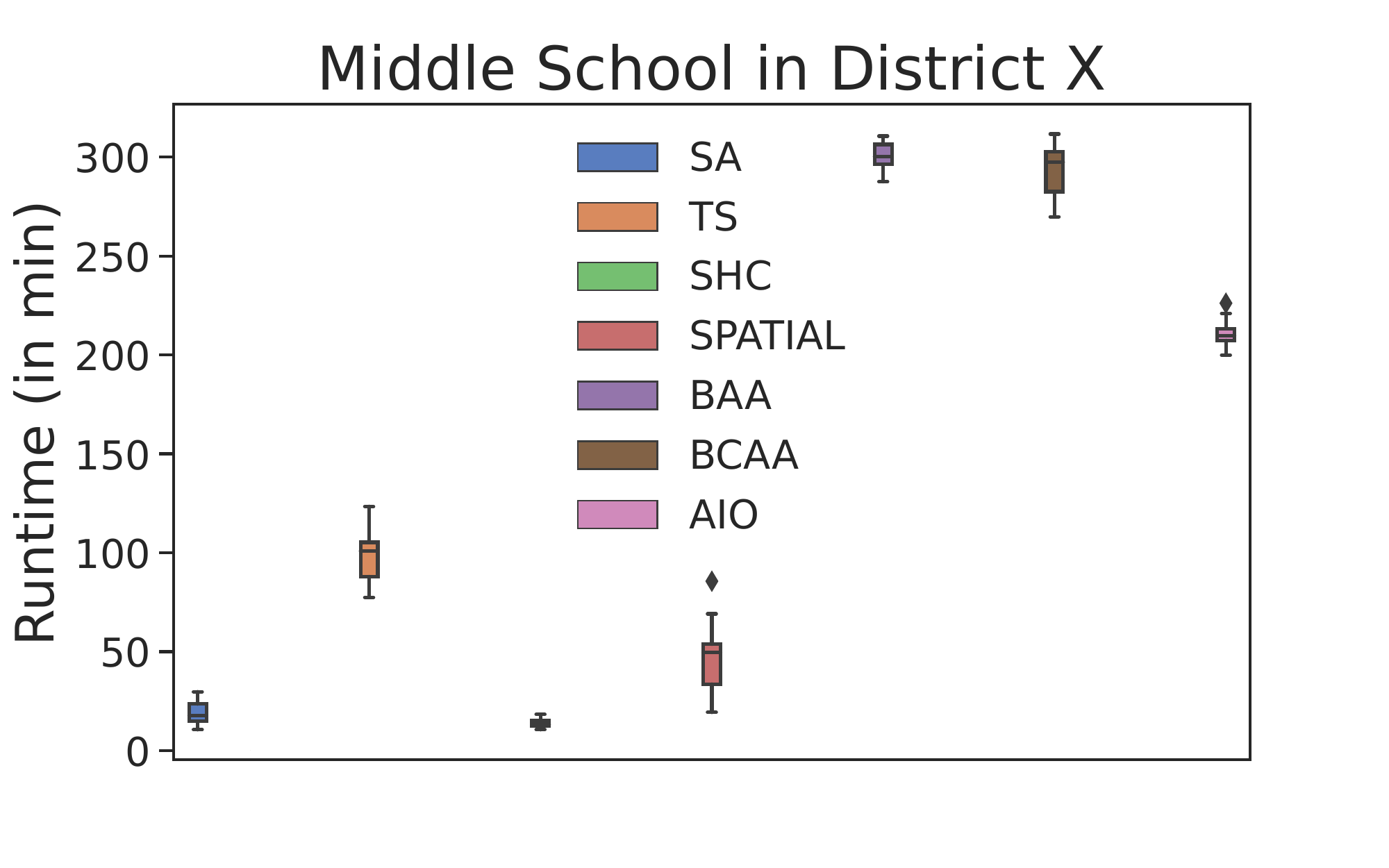}
~
        \includegraphics[width=0.495\textwidth, keepaspectratio]{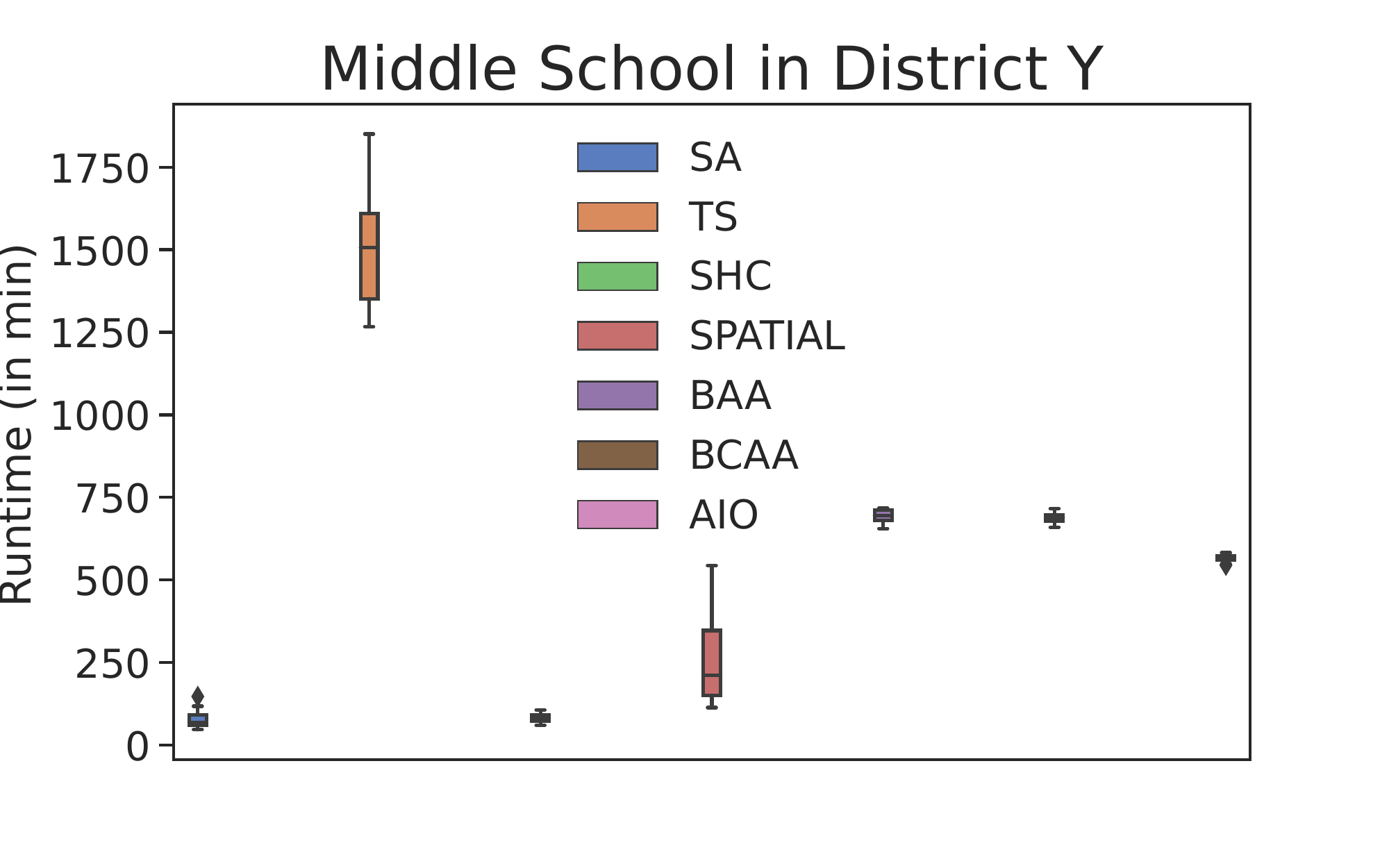}
        \caption{Middle school}
    \end{subfigure}

   \begin{subfigure}[b]{\linewidth}
        \includegraphics[width=0.495\textwidth, keepaspectratio]{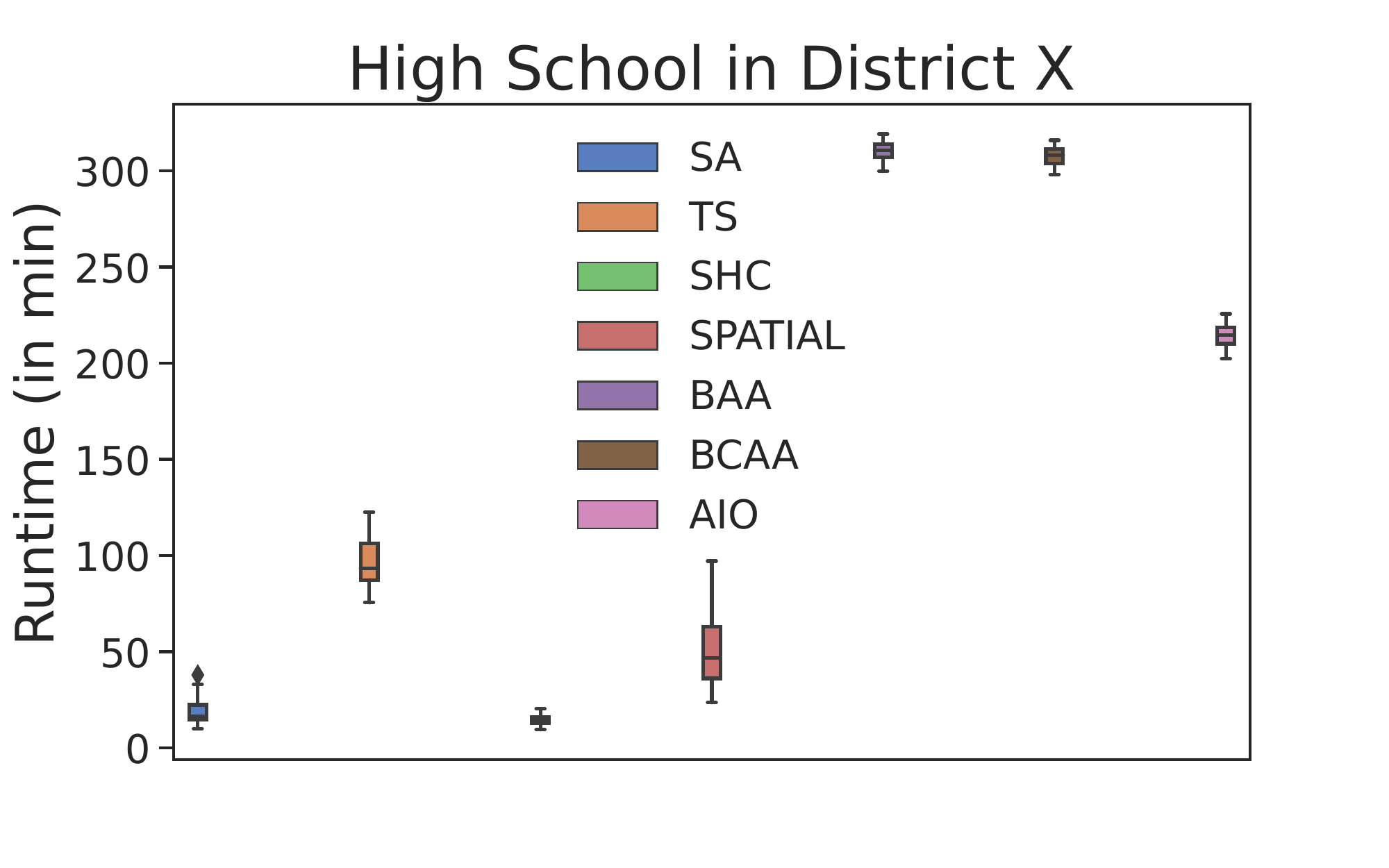}
~
        \includegraphics[width=0.495\textwidth, keepaspectratio]{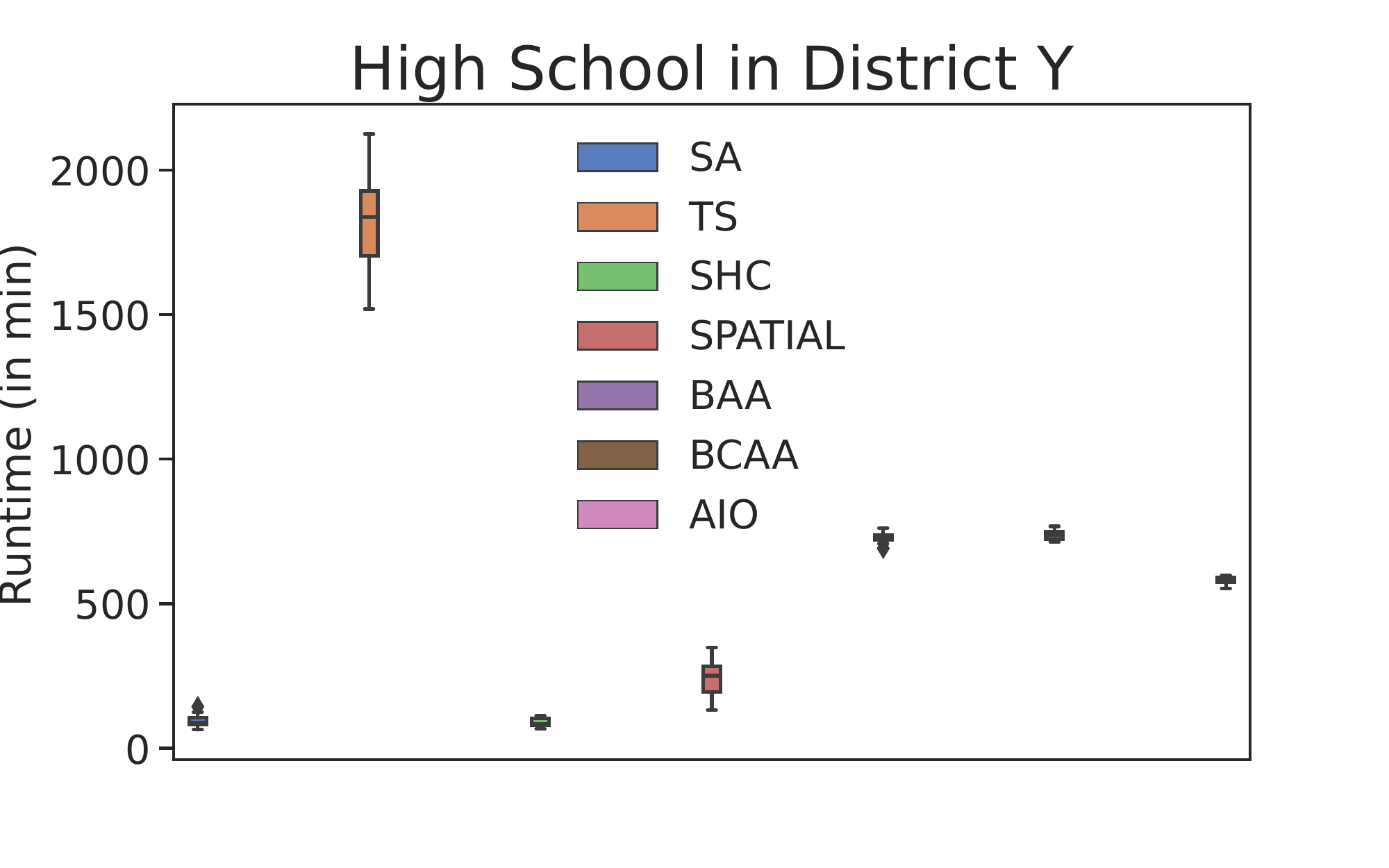}
        \caption{High school}
    \end{subfigure}
     \caption{The wall-clock time of the different methods reported for all the problem instances. The increase in the problem size causes significant increase in run time. This can be attributed to the combinatorial explosion in the search space.}
    \label{fig:runtime}
    \end{figure}

To analyze the computational complexity of these methods, we plotted the wall-clock time for all trial runs in the form of error plots in  \cref{fig:runtime}.
Run-time analysis reveals that for the smaller District~\lcps, \algo~takes longer time than the local search algorithms SA and SHC. TS is comparatively more expensive than any of the local search methods but much quicker than the sampling-based techniques that simulate 10 million flips.  Interestingly, for the larger District \fcps, TS takes the longest time and shows wide variation in the run time. This high variance can be attributed to a multitude of factors: strong dependence on initial solution, randomized order of search moves, and the tendency to get trapped in local optima. We noticed that the methods applying local search do suffer from these issues. Also, the design of local search and the computation of target compactness may be reasons for a higher computation time. Efficient algorithm design can lead to overcoming these bottlenecks.
Interestingly, the sampling-based methods are built on top of the GerryChain library which is scalable to larger-sized problem instances. 

 The calculation of the Polsby-Popper compactness measure in techniques like SA, SHC, TS and \algo, is a key computational bottleneck and it is imperative to look for faster alternatives. We posit some ways in which techniques like \algo can be made faster.
\textit{Firstly}, spatial recombination can be promoted at scale by swapping multiple nodes instead of pairwise swap in use presently. This will help to promote better exploration of the search space. 
\textit{Secondly}, finding good initial solution can help to accelerate the search significantly. For instance, replacing the random assignment of nodes (in the guided-growth phase in \Cref{sec:initialization}) with a distance-based assignment can lead to more compact territories to being with. 
\textit{Lastly}, for large-sized graphs roughly corresponding to a grid structure and having high number of nodes per subgraph, the compactness measure can be closely approximated by the edge cuts.

\begin{wrapfigure}{R}{0.5\textwidth}
      \begin{center}
        \includegraphics[keepaspectratio, width=0.5\textwidth]{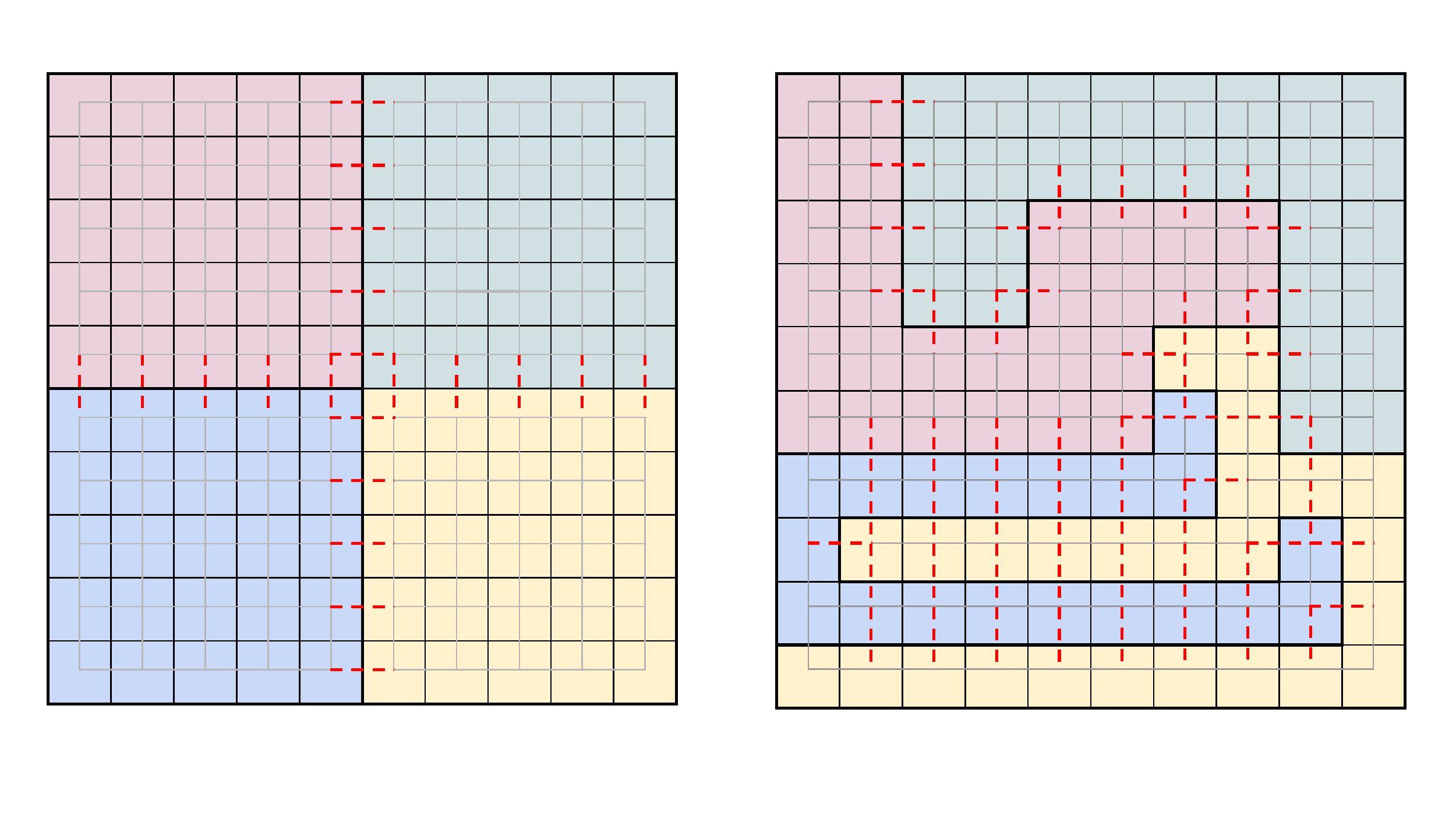}
    \caption{The dual graph of a $10 \times 10$ grid contains $180$ edges (marked in grey). The red markers indicate the edges removed to generate the corresponding partitions. Only $20$ of these edges will be removed by a partition that divides the grid into $4$ compact territories (left). However, a plan with arbitrary-shaped territories (right) could remove up to $54$ edges. This concept image is adapted from \cite{becker2021redistricting}.}
    \label{fig:cc}
    \end{center}
\end{wrapfigure}

 \Cref{fig:cc} shows that a partition with compact territories can be obtained by minimizing the number of edge cuts. The idea of edge cuts is a popular concept in graph partitioning literature where algorithms like min-cut exist~\cite{graphpartitioning}. While this approximation might not produce compact territories for problems like school redistricting, it might work very well for political districting problems~\cite{altman1998districting,becker2021redistricting}. 

As shown in  \Cref{fig:cc}, to get compact boundaries we need to minimize the number of edge cuts or retain the maximum number of existing edges while partitioning graphs. As such, the number of edges retained can be a proxy to the compactness metric used here.
Simulations on the school redistricting problem have revealed 5-10x speedup in the algorithms using local search. However, it comes at an added cost. The compactness of the boundaries suffer, this may lead to unnecessarily elongated school boundaries and consequently increased commute time for students or higher transportation costs for the school administration. However, for designing school boundaries in city blocks, most of which are organized in square blocks, the notion of edge cuts can be useful. Hence, the usage of the edge cuts for the school districting problem should be adopted with caution.


\subsection{Ablation study}
\label{sec:ablation}
\subsubsection{How effective are the search operators?
}
\label{sec:search}
To understand the effectiveness of the search operators$-$\textit{local search} and \textit{spatially-aware recombination}, we simulated 25 sample runs of~\algo~on district~\lcps~with three possible configurations by selectively activating the operators. They are as follows:
\begin{itemize}
    \item Only the local search operator is activated.
    \item Only the recombination operator is activated.
    \item Both the operators are activated.
\end{itemize}

\begin{figure}
    \centering
    \begin{subfigure}[b]{\linewidth}
        \includegraphics[width=0.495\textwidth, keepaspectratio]{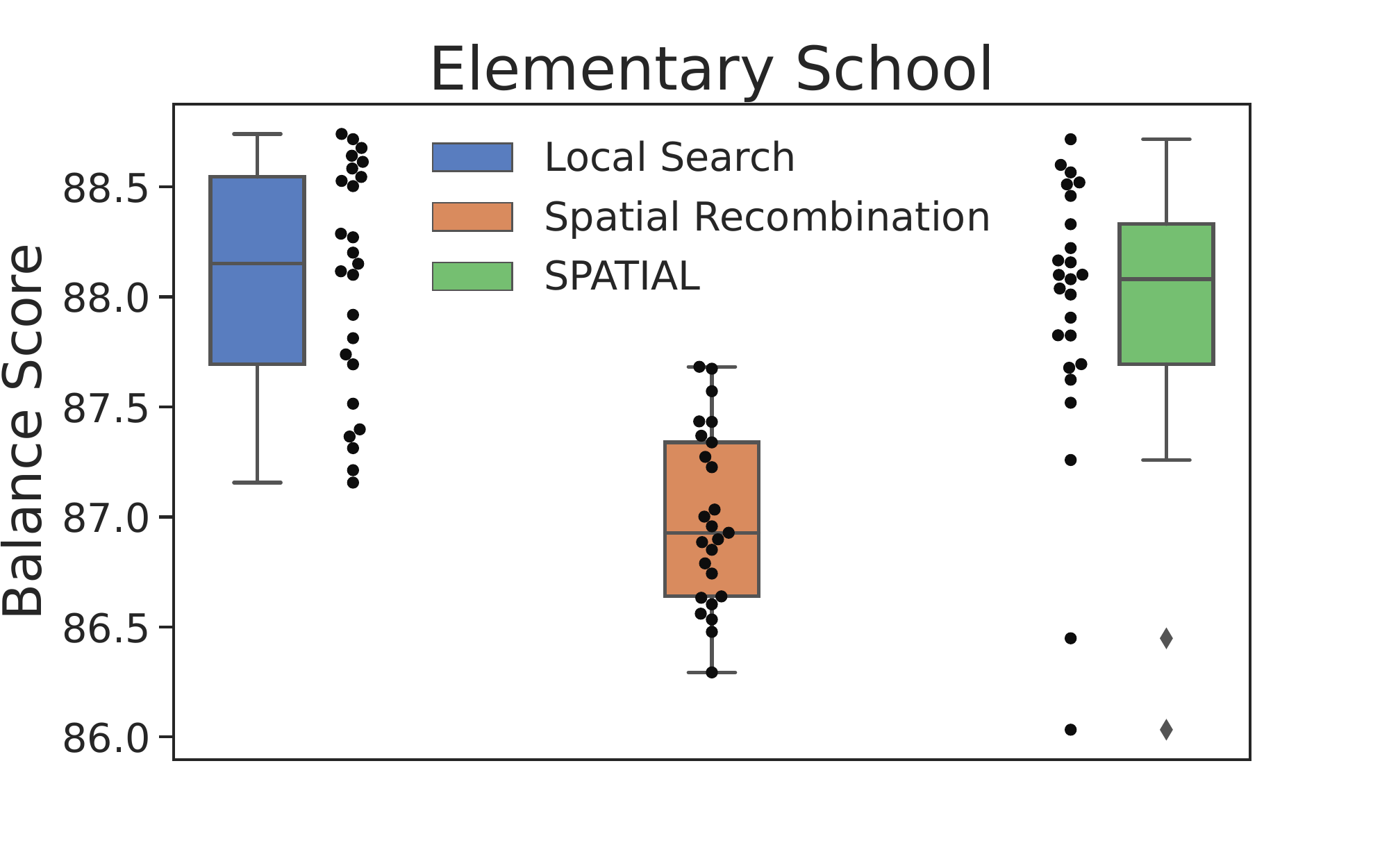}
    ~    
        \includegraphics[width=0.495\textwidth, keepaspectratio]{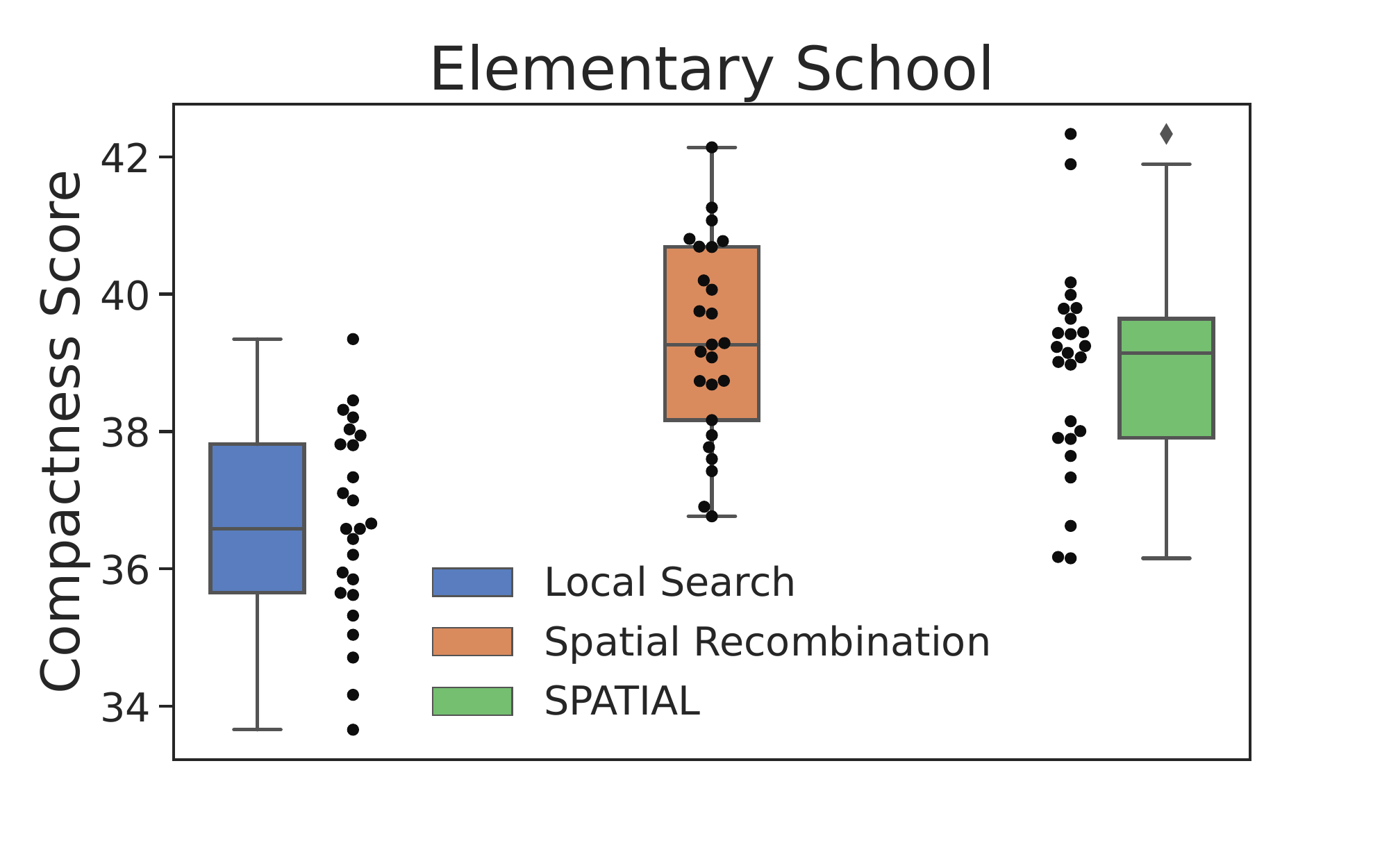}
        \caption{Elementary school}
    \end{subfigure}

   \begin{subfigure}[b]{\linewidth}
        \includegraphics[width=0.495\textwidth, keepaspectratio]{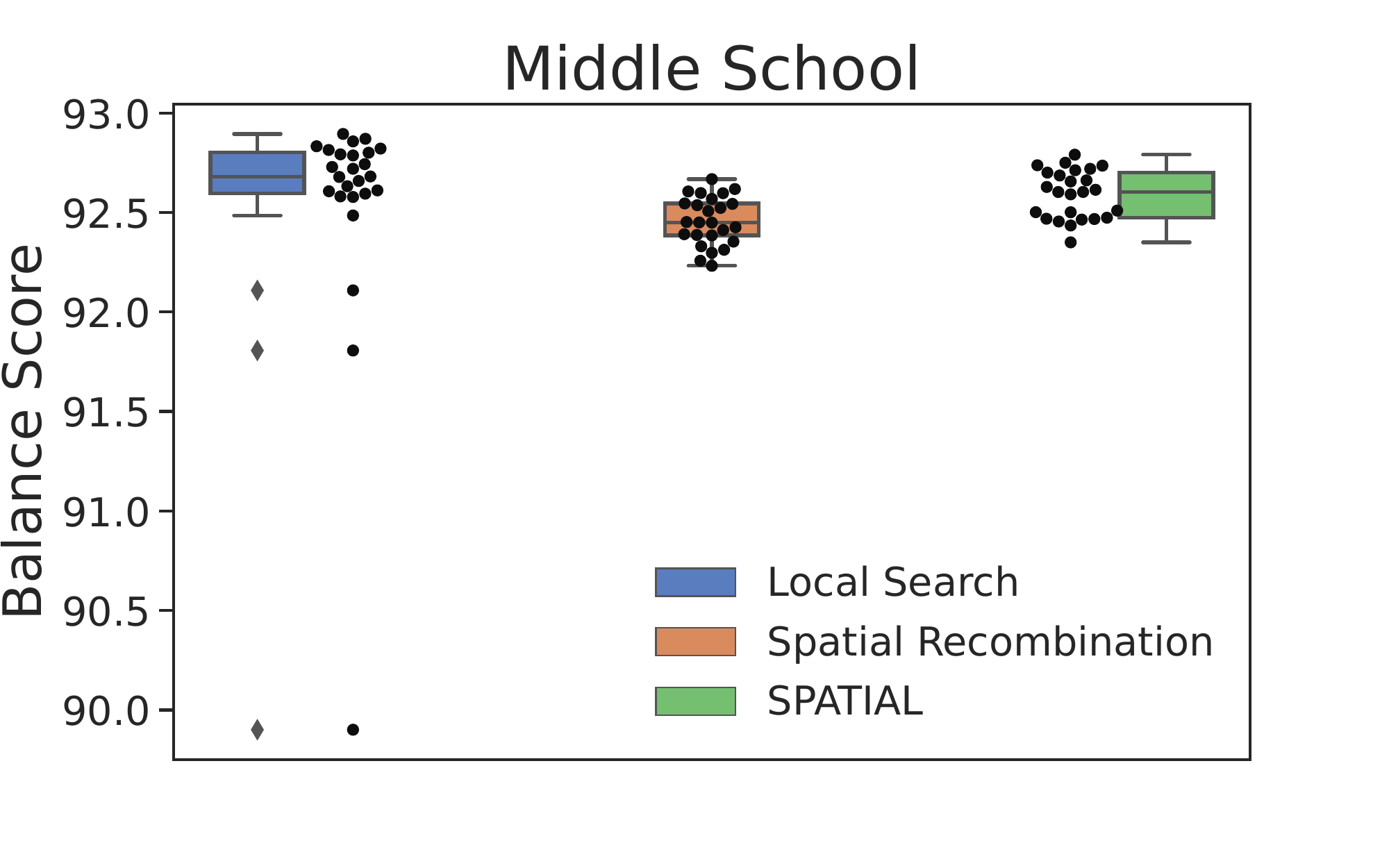}
~
        \includegraphics[width=0.495\textwidth, keepaspectratio]{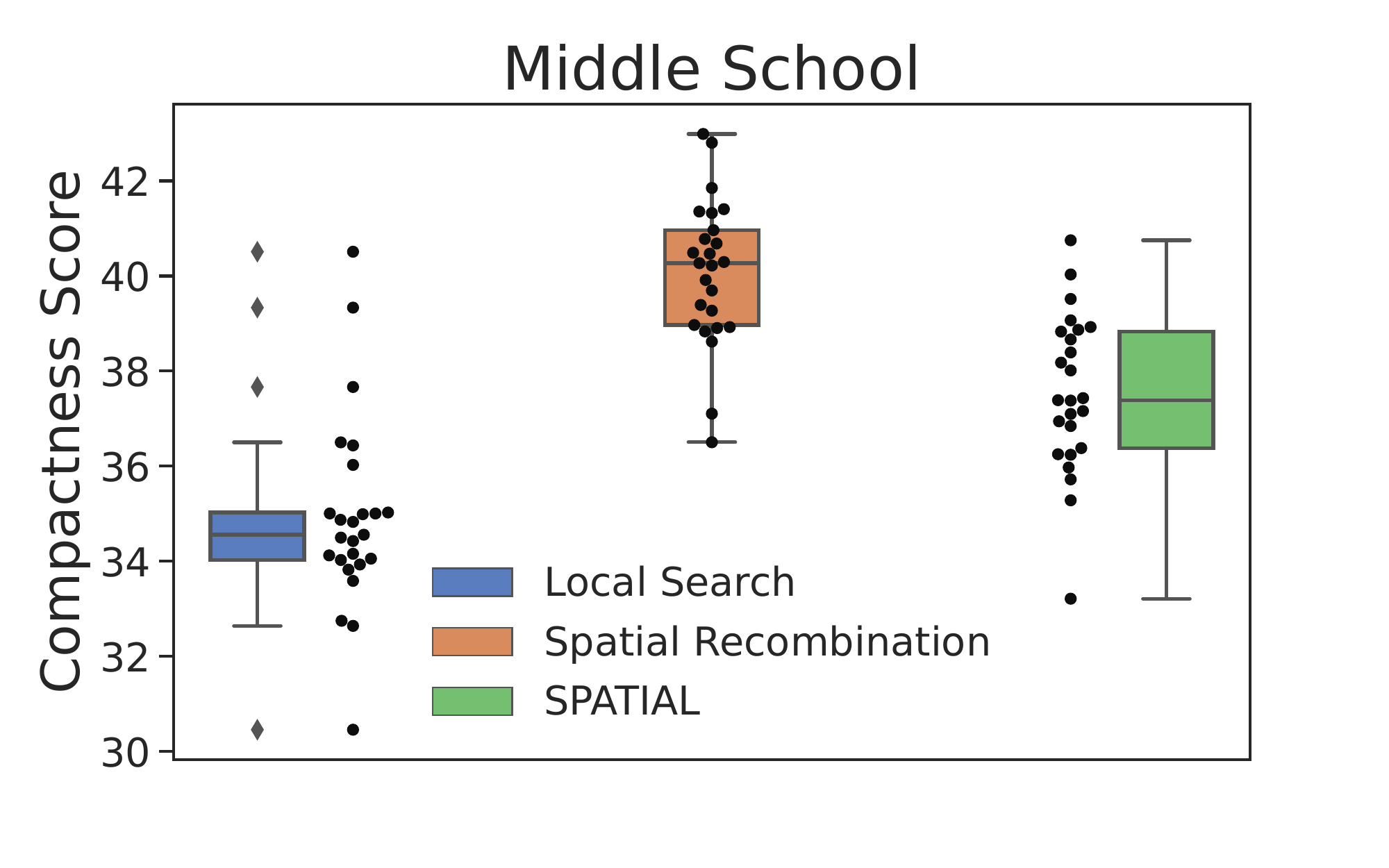}
        \caption{Middle school}
    \end{subfigure}

   \begin{subfigure}[b]{\linewidth}
        \includegraphics[width=0.495\textwidth, keepaspectratio]{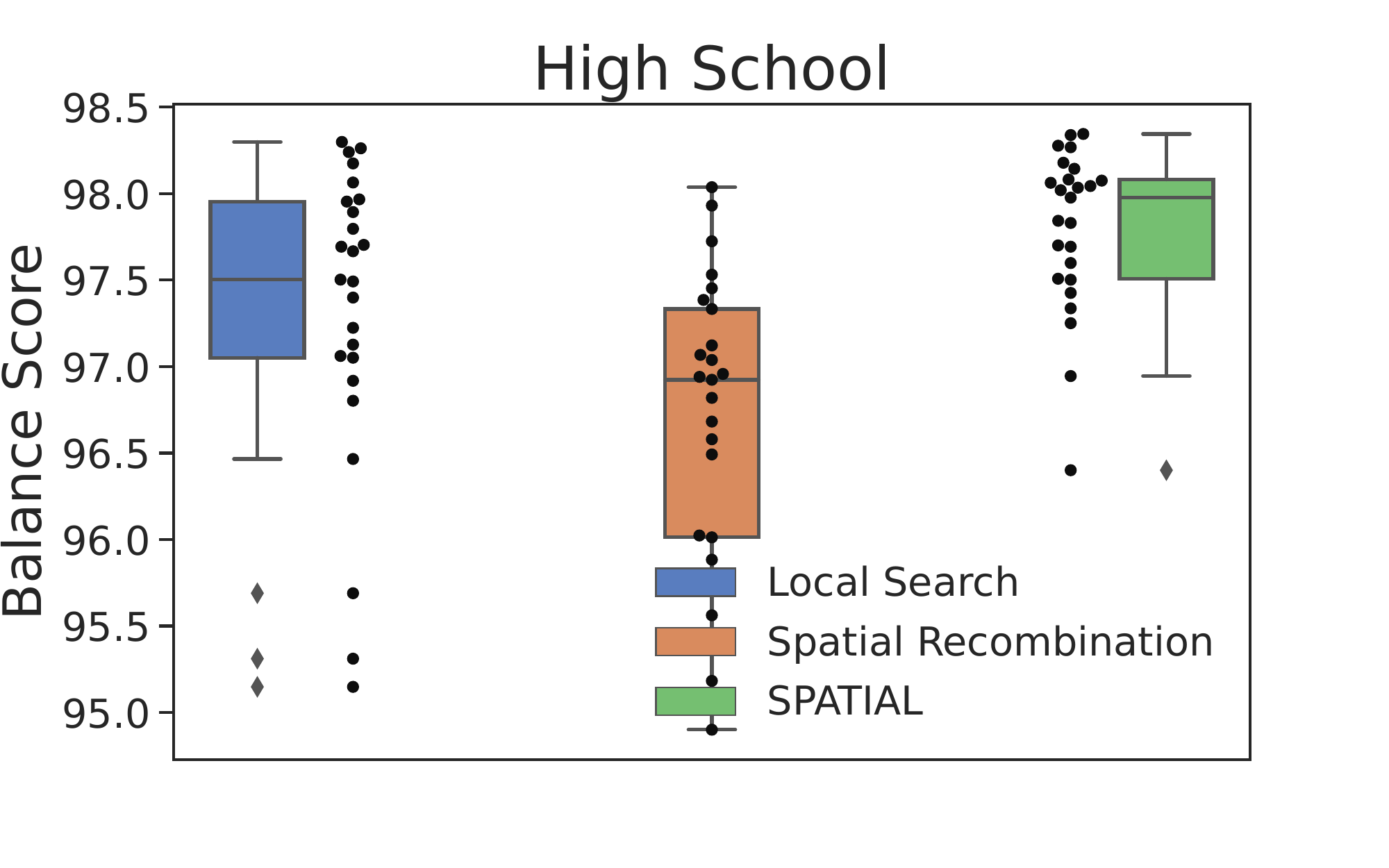}
~
        \includegraphics[width=0.495\textwidth, keepaspectratio]{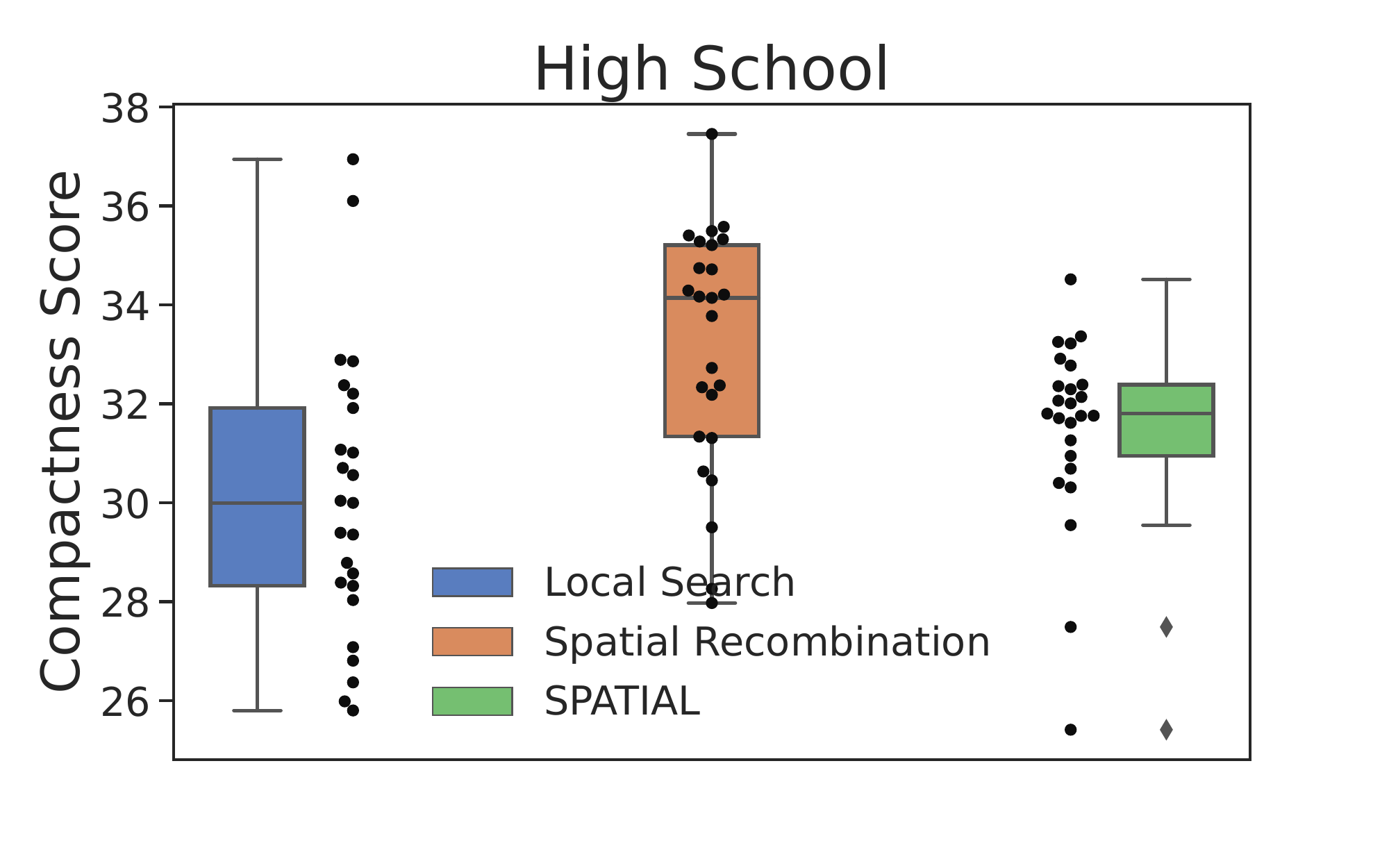}
        \caption{High school}
    \end{subfigure}
     \caption{The performance metrics obtained by the different combinations of search operators in~\algo~. We observed that the combined effect of both the operators resulted in better quality solutions.}
    \label{fig:operators}
    \end{figure}

The results of each configuration are depicted in Figure~\ref{fig:operators} as point estimates of evaluation metrics and their corresponding error plots. For fair comparison, we seeded the random numbers to ensure that the starting solutions are similar for the different configurations. Thus, any difference in the performance of these configurations can solely be attributed to the search operators.
The local search operator helps in bringing about improvement in balance scores whereas the recombination operator is responsible for high compactness scores. As mentioned earlier, recombination is less greedy than local improvement and is able to find (better) intermediate solutions beyond the immediate neighborhood of the incumbent solution. Interestingly, when both the operators are active, we noticed that the combined effect of the operators resulted in overall improvement in the quality of solutions. While local improvement resulted in exploitation of the decision space, recombination caused controlled exploration of the search space. Both are important for designing practicable school boundaries.

\subsubsection{What effect does population size have on performance?}
\label{sec:popsize}
\algo~employs a population of trial solutions for solving \spatial s. To study the effect of population size $Np$ on the performance of the algorithm, 15 runs of \algo~are simulated on the three test cases of district~\lcps~for different values of population size, $Np \in \{10, 20, 30, 40, 50\}$. 
We observed increasing the population size does not always translate to improvement in performance. This is a classic case of exploration-exploitation trade-off prevalent in optimization algorithms. A higher population size contributes to diversification (exploration) since more solutions can be distributed over the decision space. To enhance intensification (exploitation), the number of iterations need to increase in proportion to the increase in population size. This would result in longer execution time. We came across another related observation$-$ population-based metaheuristics without recombination operation may not benefit from a large population size, especially when the objective function is multi-modal and the functional landscape has attraction basins with a local minimum. This happened since the solutions present in this attraction basins would quickly eliminate other promising solutions (that could have led to a better optima) due to fitness-based replacement of the solutions. The objective function in \eqref{eq:objectivefunction} presents similar challenges as it is multi-modal in nature. The recombination operation (along with the repair operation) helps to balance such basins of attraction and preserve the solution diversity. Our observation is in line with recent findings by~\citet{chen2012large}.

\subsection{Using \algo~in real-life planning}
\label{sec:planning}
In \cref{sec:initialization}, we show how algorithms like \algo$^*$ can be used in redesigning the school boundaries in order to arrive at an alternative districting plan. Here, we analyze the practical implications of using \algo$^*$ in redistricting the school boundaries of both the districts. The evaluation entails including additional metrics beyond the ones defined in~\cref{sec:evaluation}.
\begin{itemize}
    \item \textbf{Distance-based metrics:} As an alternative to compactness, we do include some distance-based measures to get an idea of the distance traveled by students to reach schools.
        \begin{itemize}
            \item Mean distance traveled: This is the average distance a student travels to reach their respective schools. For computing this metric, we weighted the distance between the centroid of the SPA and the assigned school location with the student population in the SPA.
            \item Maximum distance traveled: This is the maximum distance between the centroid of a SPA and its assigned school.
        \end{itemize}
    
    \item \textbf{Balance-based metrics:} In addition to balance metric, we included three more metrics to highlight what proportion of schools are balanced, overcrowded or under-enrolled.
    \begin{itemize}
            \item Number of balanced schools: The number of schools in which the attending student population is between 80-120\% of the school's program capacity.
            \item Number of under-enrolled schools: The number of schools in which the attending student population is below 80\% of the school's program capacity.
            \item Number of overcrowded schools: The number of schools in which the attending student population is above 120\% of the school's program capacity.
        \end{itemize}
    
    \item \textbf{Ethical metrics:} Displacing students should be minimized since these students lose social ties to their cohorts. Thus, assessing a plan in terms of the social impact it may have on students is equally important. We show how many students may get displaced if a given plan is implemented.
\end{itemize}

\begin{table}
\centering
\caption{Comparing the existing plans with the automated plans of District~\lcps.}
\footnotesize
\begin{tabular}{@{}l|c|c|c|c|c|c@{}}
\toprule
    & \multicolumn{2}{|c|}{Elementary School}
    & \multicolumn{2}{|c|}{Middle School}
    & \multicolumn{2}{|c}{High School}
    \\
    \cline{2-7}
    & \algo$^*$
     & \textsf{Existing}
    & \algo$^*$
    & \textsf{Existing}
    & \algo$^*$
    & \textsf{Existing}
    \\
     \hline
    Compactness score
    & \textbf{42.19}
    & 32.53
    & \textbf{41.15}
    & 26.77
    & \textbf{36.07}
    & 27.35
    \\
    Mean distance traveled (in miles)
    & 0.75
    & 0.75
    & \textbf{1.24}
    & 1.28
    & 1.63
    & \textbf{1.52}
    \\
    Max distance traveled (in miles)
    & \textbf{9.84}
    & 11.75
    & \textbf{15.51}
    & 17.14
    & 15.84
    & \textbf{14.19}
    \\
    \midrule
    Balance score
    & \textbf{88.60}
    & 83.50
    & \textbf{92.78}
    & 89.74
    & \textbf{98.27}
    & 87.08
      \\
    Number of balanced schools
    & \textbf{42}/57
    & 31/57
    & \textbf{15}/17
    & 14/17
    & \textbf{16}/16
    & 14/16
      \\
    (in \%)
    & \textbf{73.68}
    & 54.39
    & \textbf{88.24}
    & 82.35
    & \textbf{100.00}
    &	87.50
    \\ 
    Number of under-enrolled schools
    & \textbf{0}/57
    & 1/57
    & 0/17
    & 0/17
    & 0/16
    & 0/16
    \\
    (in \%)
    & \textbf{0.00}
    & 1.75
    & 0.00
    & 0.00
    & 0.00
    & 0.00
    \\ 
    Number of overcrowded schools
    & \textbf{15}/57
    & 25/57
    & \textbf{2}/17
    & 3/17
    & \textbf{0}/16
    & 2/16
    \\
    (in \%)
    & \textbf{26.32}
    & 43.86
    & \textbf{11.76}
    & 17.65
    & \textbf{0.00}
    & 12.50
    \\
    \midrule
    Students displaced
    & 8253/37521
    & $-$
    & 2376/20059
    & $-$
    &3269/26728
    &$-$
    \\
    (in \%)
    & 22.00
    & $-$
    & 11.85
    & $-$
    & 12.23
    & $-$
    \\  
    \bottomrule
\end{tabular}
\normalsize
\label{tab:lcpsexistingvsautomated}
\end{table}

\begin{table}
\centering
\caption{Comparing the existing plans with the automated plans of District~\fcps.}
\footnotesize
\begin{tabular}{@{}l|c|c|c|c|c|c@{}}
\toprule
    & \multicolumn{2}{|c|}{Elementary School}
    & \multicolumn{2}{|c|}{Middle School}
    & \multicolumn{2}{|c}{High School}
    \\
    \cline{2-7}
    & \algo$^*$
     & \textsf{Existing}
    & \algo$^*$
    & \textsf{Existing}
    & \algo$^*$
    & \textsf{Existing}
    \\
     \hline
    Compactness score
    & \textbf{36.75}
    & 35.92
    & \textbf{34.28}
    & 27.71
    & \textbf{36.44}
    & 26.80
    \\
    Mean distance traveled (in miles)
    & 0.71
    & \textbf{0.68}
    & \textbf{1.68}
    & 1.77
    & 1.87
    & 1.88
    \\
    Max distance traveled (in miles)
    & 5.05
    & \textbf{4.13}
    & \textbf{11.75}
    & 12.26
    & \textbf{13.51}
    & 14.02
    \\
    \midrule
    Balance score
    & \textbf{95.81}
    & 82.38
    & \textbf{92.09}
    & 84.23
    & \textbf{91.90}
    & 86.95
      \\
    Number of balanced schools
    & \textbf{132}/138
    & 91/138
    & \textbf{22}/26
    & 17/26
    & \textbf{21}/24
    & 18/24
      \\
    (in \%)
    & \textbf{95.65}
    & 65.94
    & \textbf{84.62}
    & 65.38
    & \textbf{87.50}
    & 75.00
    \\ 
    Number of under-enrolled schools
    & \textbf{2}/138
    & 17/138	
    & \textbf{0}/26
    & 1/26
    & \textbf{3}/24
    & 6/24
    \\
    (in \%)
    & \textbf{1.45}
    & 12.32
    & \textbf{0.00}
    & 3.85
    & \textbf{12.50}
    & 25.00
    \\ 
    Number of overcrowded schools
    & \textbf{4}/138
    & 30/138
    & \textbf{4}/26
    & 8/26
    & 0/24
    & 0/24
    \\
    (in \%)
    & \textbf{2.90}
    & 21.74
    & \textbf{15.38}
    & 30.77
    & 0.00
    & 0.00
    \\
    \midrule
    Students displaced
    & 21891/100278
    & $-$
    & 6306/28647
    & $-$
    & 10749/59593
    & $-$
    \\
    (in \%)
    & 21.83
    & $-$
    & 22.01
    & $-$
    &  	18.04
    & $-$
    \\  
    \bottomrule
\end{tabular}
\normalsize
\label{tab:fcpsexistingvsautomated}
\end{table}

 For both the districts, we compare the best automated plan generated by \algo$^*$ against the existing plan and tabulate the results in \cref{tab:lcpsexistingvsautomated,tab:fcpsexistingvsautomated}.
The results in \cref{tab:lcpsexistingvsautomated} reveal that the automated plans of District~\lcps~have higher compactness values than the existing plan. However, that doesn't always translate to less distance traveled. In fact, the mean distance traveled by students increases in the automated plan is roughly same as in the existing plan. Interestingly, the automated plans were able to balance a greater proportion of schools thereby relieving the overcrowding/under-enrollment in the present schools. However, the improved balance comes at a cost. The last row reveals that if the automated plan is to be implemented, 22\%, 11.85\% and 12.23\% of the students in the elementary, middle and high schools, respectively, will be displaced with respect to the existing plan.
We noticed similar tendencies in \cref{tab:fcpsexistingvsautomated} for District \fcps. The only difference being the percentage of students displaced. District \fcps~has a high imbalance to begin with and the school administration deal with this issue by making use of modular classrooms for accommodating the extra students. However, these modular classrooms are expensive and incurs additional operational cost in the schools. The automated plan ended up changing the assignment of 15-20\% of the SPAs in order to achieve a better balance thereby resulting in high value of students displaced.

 This massive reshuffling of students is generally not encouraged and may only be done once in a span of 4-5 years. Usually, in such scenarios, the final say about which criteria to prioritize lies with the school planners. In doing so, they may consider other factors, including presence of geographic or man-made barriers, access to walk zones, socioeconomic diversity, the number of students displaced, and so on. In fact, in an ill-defined spatial problem like school districting, any automated plan cannot avoid modification. The common practice is to use arrive at the final plan by modifying a base plan or an automated plan based on subjective judgement.

\section{Conclusion \& Future Work}
\label{sec:conc}
This article proposes a metaheuristic-based approach for solving spatial optimization problems like school districting. We highlight the (computational) difficulty of using exact methods for solving such problems and motivate the need for sophisticated heuristics. To this end, the \algo~approach makes use of spatially cognizant search operators for seeking improved solutions by searching through the discrete search space characterized by spatial constraints. We illustrate two key points here. Firstly, we show how the idea of local search has theoretical underpinnings that flow from the idea of MCMC sampling on graph partitions. This even led us to compare additional sampling-based techniques for designing school boundaries. Secondly, the spatially-aware recombination along with the repair operation is instrumental in obtaining better quality solutions. An in-depth experimental investigation helped ascertain the efficacy of \algo~and related methods in solving spatial partitioning problems. We also highlight some existing drawbacks in the framework and provide pointers on ways to improve the framework.

Some possible research directions that can be undertaken in near future are as follows.
Firstly, modifying the recombination operation by incorporating multiple swaps with repair can aid in further exploration. Even techniques like ejection-chain methods can be helpful in this regard.
Secondly, developing a multi-objective version of \algo~by including multiple decision criteria. This will help to incorporate other ethical considerations, including socioeconomic diversity, equal opportunity, past displacements, and so on, into the algorithmic model.
Thirdly, integrating an sampling-based technique like MCMC with a population-based metaheuristic can help augment the search process by enhancing the diversity of the solutions.
Lastly, we can apply \algo~for solving similar \spatial s like political districting. This may require modifying the objective function and constraint-handling technique.

\begin{acks}
We are grateful to Jessica Gillis, Susan Hembach and Andreea Sistrunk for providing us with the data and helpful insights regarding the school boundary process.
\textbf{Disclaimer:} The views and conclusions contained herein are those of the authors and should not be interpreted as necessarily representing the official policies or endorsements, either expressed or implied, of any school board, NSF, or the U.S. Government.
\end{acks}

\bibliographystyle{ACM-Reference-Format}
\bibliography{references}


\end{document}